# SUMMARY OF FOUR GENERALIZED EXPONENTIAL MODELS (GEM) FOR CONTINUOUS PROBABILITY DISTRIBUTIONS

**Revised Edition (Corrected)**

◀▶

Francis J. O'Brien, Jr., Ph.D.
Naval Undersea Warfare Center, Division Newport
Undersea Warfare Combat Systems Department
1176 Howell St.
Newport, RI 02841-1708
e-mail: francis.j.obrien@navy.mil

╬

**May 5, 2014**

---

## Abstract

Four new probability models are derived which generalize the common univariate continuous distributions. Classical distributional measures are derived from Hoel, et al., *Introduction to Probability Theory,* 1971. Measures include probability density function, moments generating function, cumulative distribution function, derivatives, skewness, kurtosis, change of variable distributions, inverse distributions, log distributions. Maximum likelihood estimation technique is briefly outlined. Appendices describe applications. Errata sheet included.

# ERRATA, ADDENDA & CORRIGENDA

Summary of four generalized exponential models (GEM) for continuous probability distributions,
Francis J. O'Brien, Jr. arXiv:0801.2941v4 [math.GM]

| PAGE | NOW READS | CHANGE TO |
|---|---|---|
| 5,<br>3 lines from bottom | may | many |
| 8,<br>$h(Mo)$ for $Mo = 5 \pm 10\sqrt{2}$ | $h(mo) = \sqrt{\frac{2}{\pi}}\ e^{-1}$ | $h(mo) = \frac{1}{5e\sqrt{2\pi}}$ |
| 13,<br>6 lines from bottom | derived from, | derived from GR 8.352.6, |
| 28,<br>constraint for Model III transform,<br>$Y = a + bX$ | $y > 0$ | $y > a$ |
| 44,<br>formula for $\gamma$ | No subscripts on *x* variables in denominator term,<br>$p\sum_{i=1}^{p} x^n lnx$ | Add subscripts *i* |
| 53, Table A7 | --- | Add formula:<br>$\int \exp[\alpha\, x - \beta \exp(\gamma\, x)]dx = -\frac{\Gamma\left(\frac{\alpha}{\gamma}, \beta e^{\gamma\, x}\right)}{\gamma \beta^{\frac{\alpha}{\gamma}}}$<br>$\mathrm{Re}[\beta > 0, \gamma > 0]$ |
| 53, Table A7 | --- | Add formula: $$\int_u^v \exp[\alpha x-\beta\exp(\gamma x)]dx=\frac{\Gamma\left(\frac{\alpha}{\gamma},\beta e^{\gamma u}\right)-\Gamma\left(\frac{\alpha}{\gamma},\beta e^{\gamma v}\right)}{\gamma\beta^{\frac{\alpha}{\gamma}}}$$ $\mathrm{Re}[\beta>0,\gamma>0]$ |
| 54, Table A7 | --- | Add formula: $$\int_{-\infty}^{\infty} \exp[\alpha x-\beta\exp(\gamma x)]dx=\frac{\Gamma\left(\frac{\alpha}{\gamma}\right)}{\gamma\beta^{\frac{\alpha}{\gamma}}}$$ $\mathrm{Re}\left[\gamma,\beta,\frac{\alpha}{\gamma}>0\right]$ |
| 55, $\int_0^\infty x^m e^{-\frac{\beta}{x^n}}dx$ | --- | Add constraint: $\mathrm{Re}\,\gamma\neq 0,1,2,3,\ldots$ |
| 56, $\Gamma(-z)$ | --- | Change constraint: $\mathrm{Re}\,z\neq 0,1,2,3,\ldots$ |

# TABLE OF CONTENTS

Table 1. **Selected Univariate Densities Based on Algebraic-Exponential Functions Classified by Model Type I-IV**[1]

| DENSITY NAME | DOMAIN | PROBABILITY DENSITY FUNCTION $f(x)$ | PARAMETER RESTRICTION | GEM CLASSIFICATION (See Table 2 ▼) | PARAMETER TYPE | | |
|---|---|---|---|---|---|---|---|
| | | | | | RATE/ SCALE | SHAPE | LOCATION |
| Normal (Standardized) | $-\infty < x < \infty$ | $\frac{1}{\sqrt{2\pi}} e^{-\frac{1}{2}x^2}$ | none<br>Even function | I | | | |
| Error Function | $-\infty < x < \infty$ | $\frac{h}{\sqrt{\pi}} e^{-h^2 x^2}$ | $0 < h < \infty$<br>Even function | I | $h^2$ | | |
| Exponential | $0 \le x < \infty$ | $\lambda e^{-\lambda x}$ | $0 < \lambda < \infty$ | II | $\lambda$ | | |
| Gamma | $0 \le x < \infty$ | $\frac{\lambda^p}{\Gamma(p)} x^{p-1} e^{-\lambda x}$ | $0 < p < \infty$<br>$0 < \lambda < \infty$ | II | $\lambda$ | $p$ | |
| Weibull | $0 \le x < \infty$ | $abx^{b-1} e^{-ax^b}$ | $0 < a < \infty$<br>$0 < b < \infty$ | II | $a$ | $b$ | |
| Chi-Square | $0 \le x < \infty$ | $\frac{x^{\frac{\nu}{2}-1}}{\Gamma\left(\frac{\nu}{2}\right) 2^{\frac{\nu}{2}}} e^{-\frac{x}{2}}$ | $0 < \nu < \infty$ | II | | $\nu$ | |
| Rayleigh | $0 \le x < \infty$ | $\frac{x e^{-\frac{x^2}{2\sigma^2}}}{\sigma^2}$ | $0 < \sigma < \infty$ | II | $\sigma$ | | |
| Maxwell | $0 \le x < \infty$ | $\sqrt{\frac{2}{\pi}} x^2 e^{-\frac{1}{2}x^2}$ | none | II | | | |

[1] Most densities selected from (a): Hoel, et al. *Introduction to Probability Theory,* Houghton, Boston, 1971 and (b) p. 930 in M. Abramowitz and I.A. Stegun, *Handbook of Mathematical Functions with Formulas, Graphs, and Mathematical Tables,* Washington, D.C., Government Printing Office, 1964. A larger collection is found in Johnson, et al., 1994-1995.

| | | | | | | | |
|---|---|---|---|---|---|---|---|
| Nakagami | $0 \le x < \infty$ | $\frac{2\left(\frac{\mu}{\omega}\right)^{\mu}}{\Gamma(\mu)} x^{2\mu-1} e^{-\left(\frac{\mu}{\omega}x^2\right)}$ | $0<\mu<\infty$<br>$0<\omega<\infty$ | II | $\omega$ | $\mu$ | |
| Generalized Gamma (Stacy, 1962) | $0 \le x < \infty$ | $\frac{x^{d-1}}{a^d \Gamma\left(\frac{d}{p}\right)} e^{-\left(\frac{x}{a}\right)^p}$ | $0<d<\infty$<br>$0<a<\infty$<br>$0<p<\infty$ | II | $a$ | $p, d$ | |
| Transformed Gamma (Venter, 1983) | $0 \le x < \infty$ | $\frac{\alpha\lambda}{\Gamma(r)} (\lambda x)^{\alpha r-1} e^{-(\lambda r)^{\alpha}}$ | $0<\lambda<\infty$<br>$0<\alpha<\infty$<br>$0<r<\infty$ | II | $\lambda$ | $\alpha, r$ | |
| Normal (Non Standardized) | $-\infty < x < \infty$ | $\frac{1}{\sigma\sqrt{2\pi}} e^{-\frac{1}{2}\left(\frac{x-\mu}{\sigma}\right)^2}$ | $-\infty<\mu<\infty$<br>$0<\sigma<\infty$<br>Even function | III | $\sigma$ | | $\mu$ |
| Lognormal[2] | $0 \le x < \infty$ | $\frac{e^{-\frac{1}{2}\left(\frac{\ln x-\mu}{\sigma}\right)^2}}{x\sigma\sqrt{2\pi}}$ | $-\infty<\mu<\infty$<br>$0<\sigma<\infty$<br>Even function | III | $\sigma$ | | $\mu$ |
| Laplace | $-\infty < x < \infty$ | $\frac{1}{2b} e^{-\left\lvert\frac{x-a}{b}\right\rvert}$ | $-\infty<a<\infty$<br>$0<b<\infty$<br>Even function | III | $b$ | | $a$ |
| Exponential (2-Parameter) | $\gamma \le x < \infty$ | $\lambda e^{-\lambda(x-\gamma)}$ | $0<\lambda<\infty$<br>$0<\gamma<\infty$ | IV | $\lambda$ | | $\gamma$ |
| Gamma (3-parameter) | $\gamma \le x < \infty$ | $\frac{\lambda^p}{\Gamma(p)} (x-\gamma)^{p-1} e^{-\lambda(x-\gamma)}$ | $0<\lambda<\infty$<br>$0<p<\infty$<br>$0<\gamma<\infty$ | IV | $\lambda$ | $p$ | $\gamma$ |

[2] Lognormal is a special case. It is a log-transformed distribution based on the normal (non standardized) distribution. See Table A5. The lognormal is classified as Model III based on its distributional properties.

| | | | | | | | |
|---|---|---|---|---|---|---|---|
| Weibull (3-parameter) | $\gamma \le x < \infty$ | $ab(x-\gamma)^{b-1} e^{-a(x-\gamma)^b}$ | $0 < a < \infty$<br>$0 < b < \infty$<br>$-\infty < \gamma < \infty$ | IV | $a$ | $b$ | $\gamma$ |
| Pearson (Type III) | a $\le x < \infty$ | $\frac{1}{b\Gamma(p)}\left(\frac{x-a}{b}\right)^{p-1} e^{-\left(\frac{x-a}{b}\right)}$ | $-\infty < a < \infty$<br>$0 < b < \infty$<br>$0 < p < \infty$ | IV | $b$ | $p$ | $a$ |

NOTES:

- Common densities, based on GEM I-IV classification, conceived as integrands of Integral Equations of forms in Table 2, below, from which a normalizing constant is calculated to obtain PDF, $f(x)$. Model I is a symmetric density on the whole number line, $\pm\infty$, without a location parameter. Model III is the same with location parameter. Model II is on $0 \to \infty$ without a location parameter. Model IV is similar with a location parameter (which also serves as the lower limit of the domain).
  - Models I & II, Integral Equation, or modeling function, is, $g(x) = \alpha x^m \exp(-\beta x^n)$, on appropriate domain with $m, n$ comprising an even-function PDF in Model I on $\pm\infty$. (See Appendix on Even/Odd Function PDFs).
  - Models III & IV, linear transforms of Model I & II of form, $y = a + bx$; Integral Equation or modeling function, is,

    $g(x) = \alpha\left(\frac{x-a}{b}\right)^m \exp\left[-\beta\left(\frac{x-a}{b}\right)^n\right]$, on appropriate domain; Model III is an even-function PDF on $\pm\infty$.
    - Thus, standard Gaussian density, an even function, derived from modeling function,

      $$g(x) = \exp\left(-\frac{1}{2}x^2\right) \Rightarrow f(x) = \frac{1}{\sqrt{2\pi}}\exp\left(-\frac{1}{2}x^2\right)$$

  - Thus, conceptually we can think of a generalized GEMII PDF, $f(x) = \frac{g(x)}{\int_0^\infty g(x)dx} = \frac{x^m e^{-\beta x^n}}{\Gamma(\gamma)/n\beta^\gamma}, \gamma = \frac{m+1}{n} > 0, x > 0, g(x) > 0$

- NOTE: Scale parameters are called rate or inverse scale parameters in cases where they appear in numerators although some authors convert scale parameters to rate parameters by re-writing the PDFs, as in the exponential distribution, sometimes written in a form such as,
  - $f(x) = \lambda e^{-\lambda x}$, where $\lambda = \frac{1}{\theta}$

**Table 2. Summary of Probability Models: Density Functions**

| MODEL | DOMAIN | INTEGRAL EQUATION | NORMALIZING CONSTANT $c$ | PROBABILITY DENSITY FUNCTION PDF, $f(x)$ | PARAMETER RESTRICTION | NOTES |
|---|---|---|---|---|---|---|
| I<br>➢ Error<br>➢ Gaussian | $-\infty < x < \infty$ | $\int_{-\infty}^{\infty} \alpha x^m e^{-\beta x^n} dx$ | $\frac{n\beta^\gamma}{2\alpha\Gamma(\gamma)}$ | $\frac{n\beta^\gamma}{2\Gamma(\gamma)} x^m e^{-\beta x^n}$ | $\gamma = \frac{m+1}{n} > 0$<br>$m > -1, \alpha, \beta, n > 0$ | • even function PDF $f(-x) = f(x)$<br>• symmetric density |
| II<br>➢ Exponential<br>➢ Gamma/ Erlang<br>➢ Ch-sq.<br>➢ Rayleigh<br>➢ Weibull<br>➢ Maxwell<br>➢ Nakagami<br>➢ Inverse gamma<br>➢ Generalized gamma<br>➢ Transformed gamma | $0 \le x < \infty$ | $\int_{0}^{\infty} \alpha x^m e^{-\beta x^n} dx$ | $\frac{n\beta^\gamma}{\alpha\Gamma(\gamma)}$ | $\frac{n\beta^\gamma}{\Gamma(\gamma)} x^m e^{-\beta x^n}$ | $\gamma = \frac{m+1}{n} > 0$<br>$m > -1, \alpha, \beta, n > 0$ | |
| III<br>➢ Gaussian | $-\infty < x < \infty$ | $\int_{-\infty}^{\infty} \alpha \left(\frac{x-a}{b}\right)^m e^{-\beta\left(\frac{x-a}{b}\right)^n} dx$ | $\frac{n\beta^\gamma}{2\alpha b\Gamma(\gamma)}$ | $\frac{n\beta^\gamma}{2b\Gamma(\gamma)} \left(\frac{x-a}{b}\right)^m e^{-\beta\left(\frac{x-a}{b}\right)^n}$ | $\gamma = \frac{m+1}{n} > 0$<br>$m > -1, \alpha, \beta, b, n > 0$<br>$-\infty < a < \infty$ | • even function PDF $f(-x) = f(x)$<br>• symmetric density<br>• linear transform of Model I, $y = a + bx$ |
| IV<br>➢ Amoroso<br>➢ Pearson (Type III)<br>➢ 2-parm. Exponential<br>➢ 3-parm. Weibull | $a \le x < \infty$ | $\int_{a}^{\infty} \alpha \left(\frac{x-a}{b}\right)^m e^{-\beta\left(\frac{x-a}{b}\right)^n} dx$ | $\frac{n\beta^\gamma}{\alpha b\Gamma(\gamma)}$ | $\frac{n\beta^\gamma}{b\Gamma(\gamma)} \left(\frac{x-a}{b}\right)^m e^{-\beta\left(\frac{x-a}{b}\right)^n}$ | $\gamma = \frac{m+1}{n} > 0$<br>$m > -1, \alpha, \beta, b, n > 0$<br>$-\infty < a < \infty$ | • linear transform of Model II, $y = a + bx$ |

NOTES:

- Derivations in O'Brien, 2008. Appendix, Table A7, lists useful integral formulas.
- Models may require new names such as Generalized Exponential Model I, …, Generalized Exponential Model IV (GEMI,…, GEMIV)
- All parameters are real variables (see Appendix, "Mathematics of Odd and Even Function PDFs" for distinction of $m$,$n$ in GEM I & II)
- PDF = $c \bullet g(x) = \frac{g(x)}{\int g(x)dx}$, where $g(x)$ is the integrand of Integral Equation on proper domain and $c$ is a normalizing constant.
- Model I & III on $\pm\infty$ must compute even functions only so that $f(-x) = f(x)$, all $x$ (symmetric unimodal & bimodal densities)
- Models III & IV are based on linear transform, $y = a + bx$, of Models I & II with a location parameter $a$ and scale parameters $\beta$ & $b$ and shape parameters $m$ & $n$ whereas GEMI, II have no location parameter
- See Appendix for Model I & II PDF and CDF plots.
- PDFs may be written with rate parameters by setting scale parameter $\beta \to \theta^{-1}$.
- The four models may be converted to inverse models by reparameterizing the modeling function, $g(x)$. For example, a Model II modeling function,

$g(x) = \frac{1}{x^m}\exp\left(-\frac{\beta}{x^n}\right)$ results in the PDF , $f(x) = \frac{n\beta^z}{\Gamma(z)}\frac{1}{x^m}\exp\left(-\frac{\beta}{x^n}\right), z = \frac{m-1}{n} > 0, \mathrm{Re}\, n, \beta > 0.$ Table A7, Appendix, lists useful integrals for calculations of other distributional properties; the integrals of the form, $\int_u^v g(x)dx$, facilitate the work with proper substitution for $u,v$ and gamma function properties.

- One example of a Model IV distribution with other parameter restrictions (Amoroso) has been shown to subsume may other distributions that come under Models I, II & III by direct substitution or as limiting forms by reparameterization; see http://threeplusone.com/pubs/technote/CrooksTechNote003.pdf. Future work will attempt this unification.

**Table 3a. Summary of Models: Moments about Origin Function**

| MODEL | DOMAIN | MOMENTS FUNCTION[3] $EX^j$ | PARAMETER RESTRICTION | NOTES |
|---|---|---|---|---|
| I | $-\infty < x < \infty$ | $\dfrac{\beta^{-j/n}\Gamma\left(\gamma+\frac{j}{n}\right)}{\Gamma(\gamma)}$ | $\gamma = \dfrac{m+1}{n} > 0$<br>$m > -1, \beta, n > 0$ | • Arguments of gamma function > 0<br>• $EX^j = 0, j$ odd $(j \geq 1)$<br>• Odd/even functions algebra for symmetric density (Hoel, et al., Vol. I, p. 178) |
| II | $0 \leq x < \infty$ | $\dfrac{\beta^{-\frac{j}{n}}\Gamma\left(\gamma+\frac{j}{n}\right)}{\Gamma(\gamma)}$ | $\gamma = \dfrac{m+1}{n} > 0$<br>$m > -1, \beta, n > 0$ | • Arguments of gamma function > 0 |
| III | $-\infty < x < \infty$ | $\dfrac{j!}{\Gamma(\gamma)}\sum_{k=0}^{\frac{j-t}{2}}\dfrac{a^{j-2k}b^{2k}\beta^{-2k/n}\Gamma\left(\gamma+\frac{2k}{n}\right)}{(2k)!(j-2k)!}$<br>$t = \begin{cases} 0, & j \text{ even} \\ 1, & j \text{ odd} \end{cases}$ | $\gamma = \dfrac{m+1}{n} > 0$<br>$m > -1, \beta, b, n > 0$<br>$-\infty < a < \infty$ | • Arguments of gamma function > 0<br>• $EX^j$ computes nonzero moments for even functions only (dropping all odd-terms), for $j = 1,2k$ $(k > 0)$ by counting mechanism $\binom{j}{2k}$ in binomial expansion<br>• Odd/even functions algebra for symmetric density (Hoel, et al., Vol. I, p. 178)<br>• May be stated with Binomial Coefficient, $\binom{j}{2k}$ |
| IV | $a \leq x < \infty$ | $\dfrac{j!}{\Gamma(\gamma)}\sum_{k=0}^{j}\dfrac{a^{j-k}b^{k}\beta^{-k/n}\Gamma\left(\gamma+\frac{k}{n}\right)}{k!(j-k)!}$ | $\gamma = \dfrac{m+1}{n} > 0$<br>$m > -1, \beta, b, n > 0$<br>$-\infty < a < \infty$ | • Arguments of gamma function > 0<br>May be stated with Binomial Coefficient, $\binom{j}{k}$<br>• |

- The integrals cannot be calculated in closed form for the Moment Generating Function, $M_X(t) = Ee^{tX}$, Hoel et al., Vol. I, p. 197, and the Characteristic Function, $\varphi_X(t) = Ee^{itX}$, Hoel et al., Vol. I, p. 200.
- Calculations by integral formulas in Table A7; GEMI uses GR 3.381.11;GEMII uses 3.326.2;GEM III & IV use change of variable, binomial expansion and GR 3.338.11.

---

[3] $EX^j = \int x^j f(x)dx$

**Table 3b. Summary of Models: Moments about Mean Function**

| MODEL | DOMAIN | MOMENTS FUNCTION[4] $E(X-\mu)^j$ | PARAMETER RESTRICTION | NOTES |
|---|---|---|---|---|
| I | $-\infty < x < \infty$ | $\dfrac{\beta^{-j/n}\Gamma(\gamma + j/n)}{\Gamma(\gamma)}$ | $\gamma = \dfrac{m+1}{n} > 0$<br>$m > -1, \beta, n > 0$ | • Arguments of gamma function > 0<br>• $E(X-\mu)^j = 0,\ j$ odd $(j \ge 1)$<br>• Odd/even functions algebra for symmetric density (Hoel, et al., Vol. I, p. 178) |
| II | $0 \le x < \infty$ | $\dfrac{j!\beta^{-j/n}}{\Gamma(\gamma)}\sum_{k=0}^{j}(-1)^{j-k}\dfrac{\left[\dfrac{\Gamma\left(\gamma+\frac{1}{n}\right)}{\Gamma(\gamma)}\right]^{j-k}\Gamma\left(\gamma+\frac{k}{n}\right)}{k!(j-k)!}$ | $\gamma = \dfrac{m+1}{n} > 0$<br>$m > -1, \beta, n > 0$ | • Arguments of gamma function > 0<br>• May be stated with Binomial Coefficient, $\binom{j}{k}$ |
| III | $-\infty < x < \infty$ | $\dfrac{b^j\beta^{-j/n}\Gamma(\gamma + j/n)}{\Gamma(\gamma)}$ | $\gamma = \dfrac{m+1}{n} > 0$<br>$m > -1, n, \beta, b > 0$ | • Arguments of gamma function > 0<br>• $E(X-\mu)^j = 0,\ j$ odd $(j \ge 1)$<br>• Odd/even functions algebra for symmetric density (Hoel, et al., Vol. I, p. 178) |
| IV | $a \le x < \infty$ | $\dfrac{j!b^j\beta^{-j/n}}{\Gamma(\gamma)}\sum_{k=0}^{j}(-1)^{j-k}\dfrac{\left[\dfrac{\Gamma\left(\gamma+\frac{1}{n}\right)}{\Gamma(\gamma)}\right]^{j-k}\Gamma\left(\gamma+\frac{k}{n}\right)}{k!(j-k)!}$ | $\gamma = \dfrac{m+1}{n} > 0$<br>$m > -1, n, \beta, b > 0$ | • Arguments of gamma function > 0<br>• May be stated with Binomial Coefficient, $\binom{j}{k}$ |

➢ See Table 3a Notes.

4 $E(X-\mu)^j = \int (x-\mu)^j f(x)dx$, where $\mu$ is given in Table 3c.

## Table 3c. Summary of Models: Distribution Means, Variances and Modes

| MODEL | DOMAIN | MEAN[5] $\mu$ | VARIANCE[6] $\sigma^2$ | MODE[7] | PARAMETER RESTRICTION |
|---|---|---|---|---|---|
| I | $-\infty < x < \infty$ | 0 | $\beta^{-2/n}\dfrac{\Gamma\left(\gamma+\frac{2}{n}\right)}{\Gamma(\gamma)}$ | $\pm\left(\dfrac{m}{n\beta}\right)^{\frac{1}{n}}$ | • Arguments of gamma function > 0<br>$\gamma=\dfrac{m+1}{n}$, $m>-1, \beta, n>0$<br>• Odd/even functions algebra for symmetric density (Hoel, et al., Vol. I, p. 178) |
| II | $0 \le x < \infty$ | $\dfrac{\beta^{-\frac{1}{n}}\Gamma\left(\gamma+\frac{1}{n}\right)}{\Gamma(\gamma)}$ | $\beta^{-2/n}\left\{\dfrac{\Gamma\left(\gamma+\frac{2}{n}\right)}{\Gamma(\gamma)}-\left[\dfrac{\Gamma\left(\gamma+\frac{1}{n}\right)}{\Gamma(\gamma)}\right]^2\right\}$<br>$=\beta^{-\frac{2}{n}}\dfrac{\Gamma\left(\gamma+\frac{2}{n}\right)}{\Gamma(\gamma)}-\mu^2$ | $\left(\dfrac{m}{n\beta}\right)^{\frac{1}{n}}$ | • Arguments of gamma function > 0<br>• $\gamma=\dfrac{m+1}{n}$, $m>-1, \beta, n>0$ |

---

5 $\mu = EX = \int xf(x)dx$

6 $\sigma^2 = EX^2 - (EX)^2 = \int x^2 f(x)dx - \mu^2$

7 Mo = $\left(\dfrac{df(x)}{dx}\right)_{x=0}$ Note: Models I and III are symmetric bimodal when $m \neq 0$; e.g., GEMI Gaussian, error function, exponential all have mode, $Mo \overset{m=0}{=} 0$. The maximum ordinate $h$ of modal value, $f'(x=0)$ or $h(M0)$, may be derived by substituting $x = Mo$ into $f(x)$; e.g., GEMII, $f'(x=0) = \dfrac{n^{1-\frac{m}{n}}\beta^{\frac{1}{n}}m^{\frac{m}{n}}e^{-\frac{m}{n}}}{\Gamma\left(\frac{m+1}{n}\right)}$ $(m\neq 0)$ and $\dfrac{n\beta^{\frac{1}{n}}}{\Gamma\left(\frac{1}{n}\right)}$ $(m=0)$; for GEMI & III, make denominator $2\Gamma\left(\frac{m+1}{n}\right)$. Thus, $h(M0)$, for exponential is $\lambda$ $(m=0, n=\gamma=1, \beta=\lambda)$. For GEMIII density, $f(x) = \dfrac{1}{10\sqrt{2\pi}}\left(\dfrac{x-5}{10}\right)^2 e^{-\frac{1}{2}\left(\frac{x-5}{10}\right)^2}$, mode is:

$Mo = 5 \pm 10\sqrt{2}$, h(Mo) $= \sqrt{\dfrac{2}{\pi}}e^{-1}$. GEMIV, Pearson Type III mode is $a + b(p-1)$.

NOTE: For symmetric densities, Models I & III, the Median = $\mu$, Hoel et al., Vol. I, p. 133, derivable from CDFs, Table 5a and error function, GR 8.350.1 (GR = Gradshteyn and Ryzhik, *Table of Integrals, Series and Products*).

NOTE: the calculation of mean, var., median & mode for a log-transformed dist. (log normal) is shown in Table A5.

| | | | | | |
|---|---|---|---|---|---|
| III | -∞ < $x$ < ∞ | $a$ | $b^2 \beta^{-2/n} \dfrac{\Gamma\left(\gamma + \frac{2}{n}\right)}{\Gamma(\gamma)}$ | $a \pm b\left(\dfrac{m}{n\beta}\right)^{\frac{1}{n}}$ | • Arguments of gamma function > 0<br>• $\gamma = \dfrac{m+1}{n}, m > -1, n, \beta, b > 0$<br>$-\infty < a < \infty$<br>• Odd/even functions algebra for symmetric density (Hoel, et al., Vol. I, p. 178) |
| IV | $a \le x < \infty$ | $a + b\beta^{-1/n} \dfrac{\Gamma\left(\gamma + 1/n\right)}{\Gamma(\gamma)}$ | $b^2 \beta^{-2/n} \left\{ \dfrac{\Gamma\left(\gamma + \frac{2}{n}\right)}{\Gamma(\gamma)} - \left[ \dfrac{\Gamma\left(\gamma + \frac{1}{n}\right)}{\Gamma(\gamma)} \right]^2 \right\}$<br>$= b^2 \beta^{-\frac{2}{n}} \left[ \dfrac{\Gamma\left(\gamma + \frac{2}{n}\right)}{\Gamma(\gamma)} \right] - (\mu - a)^2$ | $a + b\left(\dfrac{m}{n\beta}\right)^{\frac{1}{n}}$ | • Arguments of gamma function > 0<br>• $m > -1, n, \beta, b > 0$<br>$-\infty < a < \infty$ |

**Table 4. Moments Functions of Common Densities Calculated from Probability Models in Table 3a**

| DENSITY NAME | GEM MODEL | MOMENTS FUNCTION $EX^j$ | MEAN $\mu$ | VARIANCE $\sigma^2$ | PARAMETER RESTRICTION (See Table 1) |
|---|---|---|---|---|---|
| Normal (Standardized) | I | $\dfrac{j!}{2^{\frac{j}{2}}\left(\frac{j}{2}\right)!}$ | 0 | 1 | • $EX^j$ computes $j$ even only<br>• $EX^j = 0$, $j$ odd ($j \geq 1$)<br>• Odd/even functions algebra for symmetric density (Hoel, et al., Vol. I, p. 178)<br>• $\Gamma\left(\frac{j+1}{2}\right) = \dfrac{\sqrt{\pi}\, j!}{2^j\left(\frac{j}{2}\right)!}$ for $j$ + 1 odd integer |
| Error Function | I | $\dfrac{h^{-j}\, j!}{2^j\left(\frac{j}{2}\right)!}$ | 0 | $\dfrac{1}{2h^2}$ | • $EX^j$ computes $j$ even only<br>• $EX^j = 0$, $j$ odd ($j \geq 1$)<br>• Odd/even functions algebra for symmetric density (Hoel, et al., Vol. I, p. 178)<br>• $\Gamma\left(\frac{j+1}{2}\right) = \dfrac{\sqrt{\pi}\, j!}{2^j\left(\frac{j}{2}\right)!}$ for $j$ + 1 odd integer |
| Exponential | II | $j\lambda^{-j}\Gamma(j)$ | $1/\lambda$ | $1/\lambda^2$ | |
| Gamma | II | $\dfrac{\lambda^{-j}\Gamma(p+j)}{\Gamma(p)}$ | $p/\lambda$ | $p/\lambda^2$ | |
| Chi-Square | II | $\left(\frac{1}{2}\right)^{-j}\dfrac{\Gamma\left(\frac{\nu}{2}+j\right)}{\Gamma(\nu/2)}$ | $\nu$ | $2\nu$ | |
| Rayleigh | II | $\sigma^j 2^{\frac{j}{2}}\Gamma\left(1+\frac{j}{2}\right)$ | $\sigma\sqrt{\frac{\pi}{2}}$ | $\dfrac{4-\pi}{2}\sigma^2$ | |
| Weibull | II | $\frac{j}{b}a^{-j/b}\Gamma\left(\frac{j}{b}\right)$ | $\frac{1}{b}a^{-1/b}\Gamma\left(\frac{1}{b}\right)$ | $a^{-2/b}\left\{\frac{2}{b}\Gamma\left(\frac{2}{b}\right)-\left[\frac{1}{b}\Gamma\left(\frac{1}{b}\right)\right]^2\right\}$ | |
| Maxwell | II | $\left(\frac{1}{2}\right)^{-j/2}\frac{\Gamma\left(\frac{j+3}{2}\right)}{\sqrt{\pi}/2}$ | $2\sqrt{\frac{2}{\pi}}$ | $\frac{3\pi-8}{\pi}$ | |
| Normal (Non Standardized) | III | $j!\sum_{k=0}^{\frac{j-t}{2}}\frac{\mu^{j-2k}\left(\frac{\sigma^2}{2}\right)^k}{k!(j-2k)!}$ $t=\begin{cases}0, j \text{ even}\\1, j \text{ odd}\end{cases}$ | μ | $\sigma^2$ | • $EX^j$ computes nonzero moments for even function only (dropping all odd terms)<br>• Odd/even functions algebra for symmetric density (Hoel, et al., Vol. I, p. 178)<br>• $\Gamma\left(\frac{2k+1}{2}\right)=\frac{\sqrt{\pi}(2k)!}{2^{2k}k!}$ for $2k+1$ odd integer |
| Laplace | III | $j!\sum_{k=0}^{\frac{j-t}{2}}\frac{a^{j-2k}b^{2k}}{(j-2k)!}$ $t=\begin{cases}0, j \text{ even}\\1, j \text{ odd}\end{cases}$ | *a* | $2b^2$ | • $EX^j$ computes nonzero moments for even function only (dropping all odd terms) |
| Pearson (Type III) | IV | $\frac{j!}{\Gamma(p)}\sum_{k=0}^{j}\frac{a^{j-k}b^k\Gamma(p+k)}{(j-k)!}$ | $a + pb$ | $pb^2$ | |

NOTE: Confirmatory central moments & modes can be calculated from Table 3b & 3c for parameters of Table 1, such as,

- GEMI Gaussian, $E(X-\mu)^j = \frac{j!\sigma^j}{2^{\frac{j}{2}}\left(\frac{j}{2}\right)!}, j \neq 1,3,5,...$
- GEMII Weibull, $E(X-\mu)^j = j!a^{-\frac{j}{b}}\sum_{k=0}^{j}\frac{(-1)^{j-k}\Gamma\left(1+\frac{1}{b}\right)^{j-k}\Gamma\left(1+\frac{k}{b}\right)}{k!(j-k)!}$
- GEMII gamma mode, $Mo = \frac{p-1}{\lambda}$
- GEMIII Gaussian mode, $Mo = \mu$

**Table 5a. Cumulative Distribution Functions for Four Models**

$$F(x) = \int_{-\infty}^{x} f(t)dt = P(X \le x)$$

| MODEL | $F(x)$ | NOTES: |
|---|---|---|
| I | $\frac{1}{2}\left[1+\frac{\gamma\left(\nu,\beta x^n\right)}{\Gamma(\nu)}\right]$ | 1. $\nu = \frac{m+1}{n}$ in arguments<br>2. Models III/IV derivations assume two limits:<br>$\left.\begin{matrix}\lim_{x\to+\infty}\left\lvert\frac{x-a}{b}\right\rvert = +\infty \\ \lim_{x\to a}\left\lvert\frac{x-a}{b}\right\rvert = 0\end{matrix}\right\}$ useful for gamma limits: $\begin{cases}\Gamma(a,0)=\gamma(a,\infty)=\Gamma(a)\\ \Gamma(a,\infty)=\gamma(a,0)=0\end{cases}$ |
| II | $\frac{\gamma\left(\nu,\beta x^n\right)}{\Gamma(\nu)}$ | |
| III | $\frac{1}{2}\left[1+\frac{\gamma\left[\nu,\beta\left(\frac{x-a}{b}\right)^n\right]}{\Gamma(\nu)}\right]$ | 3. Models I/III are symmetric densities such that $F(\mu)=\frac{1}{2}$ (i.e., the median) |
| IV | $\frac{\gamma\left[\nu,\beta\left(\frac{x-a}{b}\right)^n\right]}{\Gamma(\nu)}$ | 4. $\gamma(\bullet,\bullet)$ is the lower incomplete gamma function defined in GR 8.350.1, $\Gamma(\bullet,\bullet)$ is upper incomplete gamma function, GR 8.350.2, and $\Gamma(\bullet)$ is the complete gamma function, GR 8.310.1.<br>5. See Appendix Fig. A2 for Model I & II cdf plots and programming notes<br>6. See Appendix, "Mathematics of Odd-Even Function PDFs" for intuitive derivations of symmetric densities<br>7. Alternative expressions possible by relation $\Gamma(\alpha)=\gamma(\alpha,x)+\Gamma(\alpha,x)$ and GR 3.381.10 (Table A7, Appendix) |

- Example: Model I (by property of symmetric density, GR 3.381.10 & GR 3.381.9), by alternative expression, $F(x)=1-\frac{1}{2}\frac{\Gamma\left(\gamma,\beta x^n\right)}{\Gamma(\gamma)}$ by GR 3.381.10

- Example: Model I Gaussian, $F(x)=\frac{1}{2}\left[1+\frac{\gamma\left(\frac{1}{2},\frac{x^2}{2}\right)}{\Gamma\left(\frac{1}{2}\right)}\right]=1-\frac{1}{2}\frac{\Gamma\left(\frac{1}{2},\frac{x^2}{2}\right)}{\Gamma\left(\frac{1}{2}\right)}=\frac{1}{2}\left[1+\Phi\left(\frac{x}{\sqrt{2}}\right)\right]$, by GR 3.381.10 & GR 3.281.9, where the probability integral $\Phi(\bullet)$ is the error function $\operatorname{erf}(\bullet)$, GR 8.250.1 & 8.251.1. Set $x \to \frac{x-\mu}{\sigma}$, and GEM III Gaussian CDF results.

- ➢ In general, numerical methods are required when the error function is not applicable $\left(\gamma \neq \frac{m+1}{2}\right)$, Models I/III (see Appendix, "CDF Plots")
- Example: Model II, Weibull, $F(x) = \frac{\gamma(1, ax^b)}{\Gamma(1)} = 1 - e^{-ax^b}$
- Quantiles Q or percentiles are derivable from CDFs as $P(X \leq x) = F(x) = Q$.
  - o The Quantile for $P(X \leq x) = F(x) = Q$ is the inverse function, $x = F^{-1}(Q)$, such as Median $= x = F^{-1}\left(\frac{1}{2}\right)(= \mu, \text{symmetric densities})$.
    - ▪ Numerical methods are required for quantiles except for distributions with $\gamma = \nu = \frac{m+1}{n} = 1$ (GEM II & IV) in which case a closed form solution is possible; the *Q*th quantile (GEMII) is then $x = \left[\frac{-\ln(1-Q)}{\beta}\right]^{\frac{1}{n}} = \left[\frac{\ln\left(\frac{1}{1-Q}\right)}{\beta}\right]^{\frac{1}{n}}$, derived from,

$F(x) = \frac{\gamma(1, \beta x^n)}{\Gamma(1)} = Q \Rightarrow 1 - e^{-\beta x^n} = Q \Rightarrow x = \left[\frac{-\ln(1-Q)}{\beta}\right]^{\frac{1}{n}}$. GEMIV solution is, $Q \Rightarrow x = a + b\left[\frac{-\ln(1-Q)}{\beta}\right]^{\frac{1}{n}}$.

Example: GEMII Distributions ~ Exponential, Rayleigh, and Weibull, all have $\gamma = \nu = \frac{m+1}{n} = 1$. The half-life (median, $Q = \frac{1}{2}$) for exponential is, $\frac{-\ln\left(\frac{1}{2}\right)}{\lambda} = \frac{\ln 2}{\lambda}$.

- Derivatives of CDF parameters given in Table 5d.
  - o See similar solution for 3-parameter transformed gamma in Hogg & Klugman, p. 225, $\frac{\partial F(x)}{\partial \alpha}$.

**Table 5b. Complementary Cumulative Distribution Functions for Four Models**

$$F_c(x) = \int_x^\infty f(t)dt = 1 - F(x) = P(X > x)$$

| MODEL | $F_c(x)$ |
|---|---|
| I | $\dfrac{\Gamma(\gamma, \beta x^n)}{2\Gamma(\gamma)}$ |
| II | $\dfrac{\Gamma(\gamma, \beta x^n)}{\Gamma(\gamma)}$ |
| III | $\dfrac{\Gamma\left[\gamma, \beta\left(\frac{x-a}{b}\right)^n\right]}{2\Gamma(\gamma)}$ |
| IV | $\dfrac{\Gamma\left[\gamma, \beta\left(\frac{x-a}{b}\right)^n\right]}{\Gamma(\gamma)}$ |

NOTES:

1. $\gamma = \dfrac{m+1}{n}$ in arguments
2. Models III/IV derivations assume two limits (see Table 5a):

   $$\lim_{x \to +\infty} \left| \frac{x-a}{b} \right| = +\infty$$

   $$\lim_{x \to a} \left| \frac{x-a}{b} \right| = 0$$

3. Models I/III are symmetric densities such that $F_c(\mu) = \dfrac{1}{2} = F(x)$
   - For symmetric densities (Models I/III), $F(-x) = 1 - F(x)$, defined in Hoel, et al., Vol. I, pp. 123-124, and intuitively derived in Appendix, "Mathematics of Odd and Even Function PDFs".
4. Alternative expressions possible by relation $\Gamma(\alpha) = \gamma(\alpha, x) + \Gamma(\alpha, x)$ & GR 3.381.10.

- Example: Model I (from GR 3.381.9), $F_c(x) = 1 - F(x) = 1 - \dfrac{1}{2}\left[1 + \dfrac{\gamma(\nu, \beta x^n)}{\Gamma(\nu)}\right]$ by alternative expression of $1 - F(x)$ in Table 5a.
- Example: Model II, Weibull, $F_c(x) = \dfrac{\Gamma(1, ax^b)}{\Gamma(1)} = e^{-ax^b}$
- Example: Model IV, Weibull, $F_c(x) = \dfrac{\Gamma\left[1, a(x-\gamma)^b\right]}{\Gamma(1)} = e^{-a(x-\gamma)^b}$
- Solutions for derivatives of $F_c(x)$ are described in Table 5c Notes.

- **Table 5c. Interval Cumulative Distribution Functions for Four Models**

$$F_I(x) = \int_u^v f(t)dt = F(v) - F(u) = P(u \le X \le v)$$

| MODEL | $F_I(x)$ |
|---|---|
| I | $\dfrac{\gamma(\nu, \beta v^n) - \gamma(\nu, \beta u^n)}{2\Gamma(\nu)}$ |
| II | $\dfrac{\gamma(\nu, \beta v^n) - \gamma(\nu, \beta u^n)}{\Gamma(\nu)}$ |
| III | $\dfrac{\gamma\left[\nu, \beta\left(\frac{v-a}{b}\right)^n\right] - \gamma\left[\nu, \beta\left(\frac{u-a}{b}\right)^n\right]}{2\Gamma(\nu)}$ |
| IV | $\dfrac{\gamma\left[\nu, \beta\left(\frac{v-a}{b}\right)^n\right] - \gamma\left[\nu, \beta\left(\frac{u-a}{b}\right)^n\right]}{\Gamma(\nu)}$ |

NOTES:

1. $\nu = \dfrac{m+1}{n}$ in arguments
2. Models III/IV derivations assume two limits for parameters *u* and *v* (see Table 5a):
   $\lim_{x \to +\infty} \left| \dfrac{x-a}{b} \right| = +\infty$
   $\lim_{x \to a} \left| \dfrac{x-a}{b} \right| = 0$
3. Models I/III are symmetric densities
4. Alternative expressions possible by relation $\Gamma(\alpha) = \gamma(\alpha, x) + \Gamma(\alpha, x)$ & GR 3.381.10.

- Example: Model I, as alternative expression, $F_I(x) = \dfrac{\Gamma(\gamma, \beta u^n) - \Gamma(\gamma, \beta v^n)}{2\Gamma(\gamma)}$ by GR 3.381.10
- Example: Gaussian Model I, $F_I(x) = \dfrac{\Phi\left(\frac{v}{\sqrt{2}}\right) - \Phi\left(\frac{u}{\sqrt{2}}\right)}{2}$, where $\Phi(\bullet)$ is the error function $\mathrm{erf}(\bullet)$
- Example: Model II, Weibull, $F_I(x) = e^{-au^b} - e^{-av^b}$, with $\nu = 1, n = b, \beta = a$
- Solutions for derivatives of $F_I(x)$ are described in Table 5c Notes.

**Table 5d. Derivatives of CDFs, $F(x)$, Table 5a**

| MODEL | $\frac{dF(x)}{d\beta}$ | $\frac{dF(x)}{dn}$ | $\frac{dF(x)}{da}$ | $\frac{dF(x)}{db}$ |
|---|---|---|---|---|
| I | $\frac{(\beta x^n)^{\gamma-1}}{2\Gamma(\gamma)} e^{-\beta x^n} \frac{d(\beta x^n)}{d\beta} = \frac{\beta^{\gamma-1} x^{m+1} e^{-\beta x^n}}{2\Gamma(\gamma)}$ | $\frac{(\beta x^n)^{\gamma-1}}{2\Gamma(\gamma)} e^{-\beta x^n} \frac{d(\beta x^n)}{dn} = \frac{\beta^{\gamma} x^{m+1} e^{-\beta x^n} \ln x}{2\Gamma(\gamma)}$ | — | — |
| II | $\frac{(\beta x^n)^{\gamma-1}}{\Gamma(\gamma)} e^{-\beta x^n} \frac{d(\beta x^n)}{d\beta} = \frac{\beta^{\gamma-1} x^{m+1} e^{-\beta x^n}}{\Gamma(\gamma)}$ | $\frac{(\beta x^n)^{\gamma-1}}{\Gamma(\gamma)} e^{-\beta x^n} \frac{d(\beta x^n)}{dn} = \frac{\beta^{\gamma} x^{m+1} e^{-\beta x^n} \ln x}{\Gamma(\gamma)}$ | — | — |
| III | $\frac{\left[\beta\left(\frac{x-a}{b}\right)^n\right]^{\gamma-1}}{2\Gamma(\gamma)} e^{-\beta\left(\frac{x-a}{b}\right)^n} \frac{d\left[\beta\left(\frac{x-a}{b}\right)^n\right]}{d\beta} = \frac{\beta^{\gamma-1}\left(\frac{x-a}{b}\right)^{m+1}}{2\Gamma(\gamma)} e^{-\beta\left(\frac{x-a}{b}\right)^n}$ | $\frac{\left[\beta\left(\frac{x-a}{b}\right)^n\right]^{\gamma-1}}{2\Gamma(\gamma)} e^{-\beta\left(\frac{x-a}{b}\right)^n}$ X $\frac{d\left[\beta\left(\frac{x-a}{b}\right)^n\right]}{dn} = \frac{\beta^{\gamma}\left(\frac{x-a}{b}\right)^{m+1}}{2\Gamma(\gamma)} e^{-\beta\left(\frac{x-a}{b}\right)^n} \ln\left(\frac{x-a}{b}\right)$ | $\frac{\left[\beta\left(\frac{x-a}{b}\right)^n\right]^{\gamma-1}}{2b\Gamma(\gamma)} e^{-\beta\left(\frac{x-a}{b}\right)^n}$ X $\frac{d\left[\beta\left(\frac{x-a}{b}\right)^n\right]}{da} = -\frac{n\beta^{\gamma}\left(\frac{x-a}{b}\right)^{m}}{2b\Gamma(\gamma)} e^{-\beta\left(\frac{x-a}{b}\right)^n}$ | $\frac{\left[\beta\left(\frac{x-a}{b}\right)^n\right]^{\gamma-1}}{2b\Gamma(\gamma)} e^{-\beta\left(\frac{x-a}{b}\right)^n}$ X $\frac{d\left[\beta\left(\frac{x-a}{b}\right)^n\right]}{db} = -\frac{n\beta^{\gamma}\left(\frac{x-a}{b}\right)^{m+1}}{2b\Gamma(\gamma)} e^{-\beta\left(\frac{x-a}{b}\right)^n}$ |
| IV | $\frac{\left[\beta\left(\frac{x-a}{b}\right)^n\right]^{\gamma-1}}{\Gamma(\gamma)}e^{-\beta\left(\frac{x-a}{b}\right)^n}$ X $\frac{d\left[\beta\left(\frac{x-a}{b}\right)^n\right]}{d\beta}=$ $\frac{\beta^{\gamma-1}\left(\frac{x-a}{b}\right)^{m+1}}{\Gamma(\gamma)}e^{-\beta\left(\frac{x-a}{b}\right)^n}$ | $\frac{\left[\beta\left(\frac{x-a}{b}\right)^n\right]^{\gamma-1}}{\Gamma(\gamma)}e^{-\beta\left(\frac{x-a}{b}\right)^n}$ X $\frac{d\left[\beta\left(\frac{x-a}{b}\right)^n\right]}{dn}=$ $\frac{\beta^{\gamma}\left(\frac{x-a}{b}\right)^{m+1}}{\Gamma(\gamma)}e^{-\beta\left(\frac{x-a}{b}\right)^n}\ln\left(\frac{x-a}{b}\right)$ | $\frac{\left[\beta\left(\frac{x-a}{b}\right)^n\right]^{\gamma-1}}{b\Gamma(\gamma)}e^{-\beta\left(\frac{x-a}{b}\right)^n}$ X $\frac{d\left[\beta\left(\frac{x-a}{b}\right)^n\right]}{da}=$ $-\frac{n\beta^{\gamma}\left(\frac{x-a}{b}\right)^{m}}{b\Gamma(\gamma)}e^{-\beta\left(\frac{x-a}{b}\right)^n}$ | $\frac{\left[\beta\left(\frac{x-a}{b}\right)^n\right]^{\gamma-1}}{b\Gamma(\gamma)}e^{-\beta\left(\frac{x-a}{b}\right)}$ X $\frac{d\left[\beta\left(\frac{x-a}{b}\right)^n\right]}{db}=$ $-\frac{n\beta^{\gamma}\left(\frac{x-a}{b}\right)^{m+1}}{b\Gamma(\gamma)}e^{-\beta\left(\frac{x-a}{b}\right)^n}$ |

NOTES:

- $\gamma=\frac{m+1}{n}$ in arguments
- $\frac{dF(x)}{dx}=F'(x)=f(x)$
- $\frac{d\gamma(a,x)}{dx}=-\frac{d\Gamma(a,x)}{dx}=x^{a-1}e^{-x}$; composite functions, $\frac{d\gamma(a,u)}{dx}=-\frac{d\Gamma(a,u)}{dx}=\frac{d\gamma(a,u)}{du}\frac{du}{dx}$ (GR 8.356.4)
- $\frac{d(\beta x)^n}{dn}=\frac{d(\beta^n)(x^n)}{dn}=\frac{d(e^{n\ln\beta})(e^{n\ln x})}{dn}=\frac{de^{n\ln(\beta x)}}{dn}=(\beta x)^n\ln(\beta x)$
  $\frac{dx^n}{dn}=x^n\ln(x)$
- Model II: from Table 5a, $\frac{dF(x)}{d\nu}=\frac{d}{d\nu}\frac{\gamma(\nu,\beta x^n)}{\Gamma(\nu)}$. First, solve for the primitive form $\frac{d}{da}\frac{\gamma(a,x)}{\Gamma(a)}$ and set $a\to\nu, x\to\beta x^n$. Also, state $\gamma(a,x)$ as the integral representation of the lower incomplete gamma function; i.e.,
  $\gamma(a,x)=x^a\int_0^1 t^{a-1}e^{-xt}dt,\quad$ real $a,x>0$. Then, $\frac{d}{da}\frac{\gamma(a,x)}{\Gamma(a)}=\frac{d}{da}\frac{1}{\Gamma(a)}x^a\int_0^1 t^{a-1}e^{-xt}dt$. Using the quotient rule and product rule of derivatives, and differentiating under the integral,

$$\frac{d}{da}\frac{\gamma(a,x)}{\Gamma(a)}=\frac{d}{da}\frac{1}{\Gamma(a)}x^a\int_0^1 t^{a-1}e^{-xt}dt=\frac{\gamma(a,x)}{x^a}\left[\frac{\Gamma(a)\frac{d}{da}x^a-x^a\frac{d}{da}\Gamma(a)}{[\Gamma(a)]^2}\right]+\frac{x^a}{\Gamma(a)}\frac{d}{da}\int_0^1 e^{(a-1)\ln t}e^{-xt}dt$$

$$= \frac{\gamma(a,x)}{x^a}\left[\frac{\Gamma(a)x^a \ln x - x^a\Gamma(a)\Psi(a)}{[\Gamma(a)]^2}\right] + \int_0^1 t^{a-1}e^{-xt}\ln t dt = \frac{\gamma(a,x)}{\Gamma(a)}[\ln x - \Psi(a)] + \int_0^1 t^{a-1}e^{-xt}\ln t dt$$

- Setting $a \to \nu, x \to \beta x^n$ is the solution for Model II CDF.
- Summarizing formulas are given in O'Brien ("Math. Summary", 2008)

- Model I derivative for primitive form is $\frac{d}{da}\frac{1}{2}\left[1+\frac{\gamma(a,x)}{\Gamma(a)}\right] = \frac{1}{2}\frac{d}{da}\frac{\gamma(a,x)}{\Gamma(a)}$ which is solved as above for Model II with $a \to \nu, x \to \beta x^n$ substituted.
  - Models III and IV derivatives are similar.
- Example: The 3-parameter transformed gamma distribution, $f(x) = \frac{\lambda^{\tau\alpha}\tau x^{\tau\alpha-1}e^{-(\lambda x)^\tau}}{\Gamma(\alpha)}$, Model II, Hogg & Klugman, p. 225;

$$F(x) = \frac{\gamma\left[\alpha,(\lambda x)^\tau\right]}{\Gamma(\alpha)}$$

$$\frac{dF(x)}{d\lambda} = \frac{\tau x(\lambda x)^{\alpha\tau-1}}{\Gamma(\alpha)}e^{-(\lambda x)^\tau}$$

- NOTE: the derivatives of the complementary CDFs, $F_c(x)$, Table 5b, are performed in a similar manner using the integral representation of the upper incomplete gamma function, $\Gamma(a,x) = x^a\int_1^\infty t^{a-1}e^{-xt}dt$ (see Table A7, Appendix). For example, the Model II derivative in primitive form is, $\frac{d}{da}\frac{\Gamma(a,x)}{\Gamma(a)} = \frac{d}{da}\frac{1}{\Gamma(a)}x^a\int_1^\infty t^{a-1}e^{-xt}dt = \frac{\Gamma(a,x)}{\Gamma(a)}[\ln x - \Psi(a)] + \int_1^\infty t^{a-1}e^{-xt}\ln t dt$.

➢ Derivatives of the interval CDFs, $F_I(x)$, Table 5c, are similar as differences of lower incomplete gamma functions.

**Table 6. Coefficients of Skewness and Kurtosis for Four Models**[8]

| MODEL | SKEWNESS | KURTOSIS | PARAMETER RESTRICTION |
|---|---|---|---|
| I | 0 | $\dfrac{\Gamma(\gamma)\Gamma\left(\gamma+\frac{4}{n}\right)}{\left[\Gamma\left(\gamma+\frac{2}{n}\right)\right]^2}$ | gamma arguments > 0<br>$\gamma=\dfrac{m+1}{n}>0\ m>-1, n>0$<br>even function |
| II | $\dfrac{\dfrac{\Gamma\left(\gamma+\frac{3}{n}\right)}{\Gamma(\gamma)}-3\dfrac{\Gamma\left(\gamma+\frac{2}{n}\right)\Gamma\left(\gamma+\frac{1}{n}\right)}{[\Gamma(\gamma)]^2}+2\left[\dfrac{\Gamma\left(\gamma+\frac{1}{n}\right)}{\Gamma(\gamma)}\right]^3}{\left\{\dfrac{\Gamma\left(\gamma+\frac{2}{n}\right)}{\Gamma(\gamma)}-\left[\dfrac{\Gamma\left(\gamma+\frac{1}{n}\right)}{\Gamma(\gamma)}\right]^2\right\}^{3/2}}$ | $\dfrac{\left[\dfrac{\Gamma\left(\gamma+\frac{4}{n}\right)}{\Gamma(\gamma)}\right]-4\dfrac{\Gamma\left(\gamma+\frac{3}{n}\right)}{\Gamma(\gamma)}\left[\dfrac{\Gamma\left(\gamma+\frac{1}{n}\right)}{\Gamma(\gamma)}\right]+6\dfrac{\Gamma\left(\gamma+\frac{2}{n}\right)}{\Gamma(\gamma)}\left[\dfrac{\Gamma\left(\gamma+\frac{1}{n}\right)}{\Gamma(\gamma)}\right]^2-3\left[\dfrac{\Gamma\left(\gamma+\frac{1}{n}\right)}{\Gamma(\gamma)}\right]^4}{\left\{\dfrac{\Gamma\left(\gamma+\frac{2}{n}\right)}{\Gamma(\gamma)}-\left[\dfrac{\Gamma\left(\gamma+\frac{1}{n}\right)}{\Gamma(\gamma)}\right]^2\right\}^2}$ | gamma arguments > 0<br>$\gamma=\dfrac{m+1}{n}>0$<br>$m>-1, n>0$ |
| III | SAME AS MODEL I | SAME AS MODEL I | SAME AS MODEL I<br>even function |
| IV | SAME AS MODEL II | SAME AS MODEL II | SAME AS MODEL II |

- NOTE: derived from moments relations in Tables 3a and 3b.
- Example: Kurtosis for Model I with $\gamma=\frac{1}{2}, n=2$, is, $\dfrac{\Gamma\left(\frac{1}{2}\right)\Gamma\left(\frac{5}{2}\right)}{\left[\Gamma\left(\frac{3}{2}\right)\right]^2}=3$ (Gaussian)

---

[7] Definitions. Skewness: $\gamma_1=\dfrac{E(X-\mu)^3}{\sigma^3}$. Kurtosis: $\gamma_2=\dfrac{E(X-\mu)^4}{\sigma^4}$.

# APPENDIX

**Figure A1.**

# PDF Plots

## MODEL I
### (Even Function, Symmetric Density)

*f*(*x*)

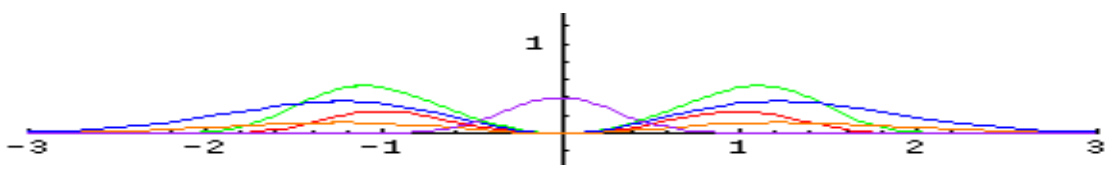

$\mu = 0$

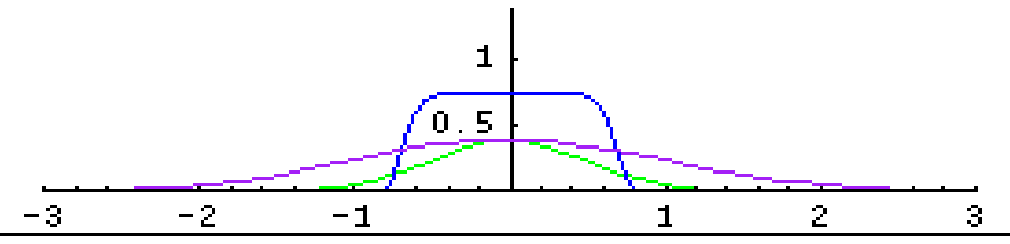

▲Plot Parameters

Top Curves

| Distribution | $m$ | $n$ | $\beta$ |
|---|---|---|---|
| Red | 2 | 2 | ½ |
| Green | 2 | 4 | 1/4 |
| Blue | 2 | 2/3 | 2 |
| Purple (Gaussian) | 0 | 2 | ½ |
| Orange | .01 | 2 | 3 |

Bottom Curves (unimodal)

| | | | |
|---|---|---|---|
| Green | 0 | 2 | 2 |
| Purple (Gaussian) | 0 | 2 | ½ |
| Blue | 0 | 10 | 33 |

$$f(x) = \frac{n\beta^{\gamma}}{2\Gamma(\gamma)} x^{m} e^{-\beta x^{n}}$$

$$\gamma = \frac{m+1}{n} > 0$$

$$m > -1, \beta, n > 0$$

## MODEL II

*f*(*x*)

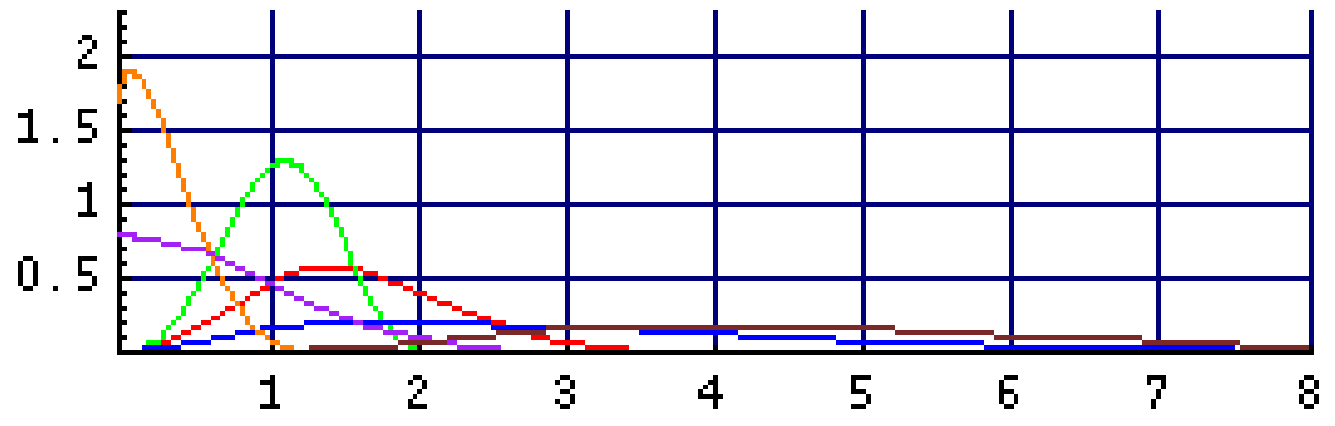

▲Plot Parameters

| Distribution | $m$ | $n$ | $\beta$ |
|---|---|---|---|
| Red | 2 | 2 | ½ |
| Green | 2 | 4 | 1/4 |
| Blue | 2 | 2/3 | 2 |
| Purple | 0 | 2 | ½ |
| Orange | .01 | 2 | 3 |
| Brown | 4 | ½ | 1 |

▲Create other distributions with high left (negative) skewness **TBD**

$$f(x) = \frac{n\beta^{\gamma}}{\Gamma(\gamma)} x^{m} e^{-\beta x^{n}}$$

$$\gamma = \frac{m+1}{n} > 0$$

$$m > -1, \beta, n > 0$$

**Figure A2.**

# CDF Plots & Numerical Methods

| MODEL I (Even Function, Symmetric Density) | MODEL II |
|---|---|
| F(x) | F(x) |

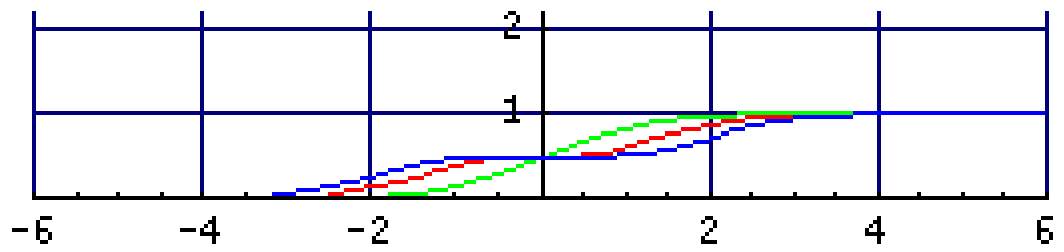

$\mu = 0$

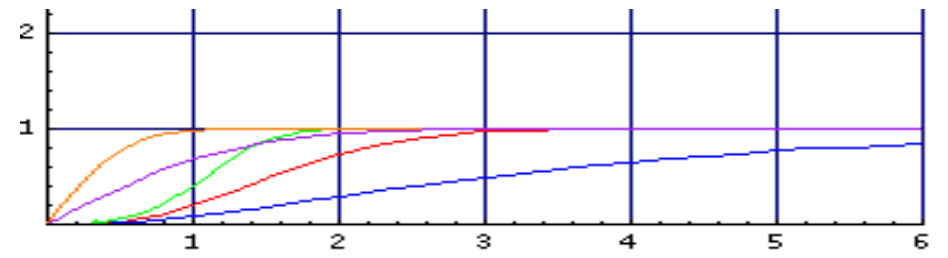

**MODEL I** ▲Plot Parameters

| Distribution | $m$ | $n$ | $\beta$ |
|---|---|---|---|
| Red | 2 | 2 | ½ |
| Green (Gaussian) | 0 | 2 | ½ |
| Blue | 4 | 2 | ½ |

**MODEL II** ▲Plot Parameters
Same as PDF (excluding Brown)

- $0 \le F(x) \le 1$.
- $F(x)$ also called Failure Distribution or Unreliability Distribution
- NOTE: All Model I curves are non-decreasing and cross at the zero point $\left(F(\mu) = F(0) = \frac{1}{2}\right)$, as expected
    - o Bimodal curves (Red & Blue) go flat near $x = 0$, as expected.
- NOTE: Model I CDF functions could not be plotted in form, $F(x) = 1 - \frac{1}{2}\frac{\Gamma(\gamma, \beta x^n)}{\Gamma(\gamma)}, \gamma = \frac{m+1}{2}$.
    - The equivalent functions in terms of error function were required, from derived equations (O'Brien, 2007, Sect. 8.359)—

$$\Gamma\left(\frac{1}{2}, x\right) = \sqrt{\pi} - \sqrt{\pi}\,\mathrm{erf}\left(\sqrt{x}\right)$$

$$\Gamma\left(\frac{3}{2}, x\right) = \frac{1}{2}\Gamma\left(\frac{1}{2}, x\right) + e^{-x}x^{\frac{1}{2}}$$

$$\Gamma\left(\frac{5}{2}, x\right) = \frac{3}{4}\Gamma\left(\frac{1}{2}, x\right) + e^{-x}x^{\frac{1}{2}}\left(\frac{3}{2} + x\right)$$

$$\Gamma\left(\frac{7}{2}, x\right) = \frac{15}{8}\Gamma\left(\frac{1}{2}, x\right) + e^{-x}x^{\frac{1}{2}}\left(\frac{15}{4} + \frac{5}{2}x + x^2\right)$$

$$\vdots$$

$$\Gamma\left(\frac{n}{2}, x\right) = \frac{(n-2)!!}{2^{\frac{n-3}{2}}}\left[\frac{1}{2}\Gamma\left(\frac{1}{2}, x\right) + e^{-x}x^{\frac{1}{2}}\sum_{s=0}^{\frac{n-3}{2}}\frac{(2x)^s}{(2s+1)(2s-1)!!}\right] \qquad [n = 1,3,5,\ldots]$$

where, $(n-2)!! = 1\cdot3\cdot5\ldots(n-2); (-1)!! = 1$.

- Error Function properties, useful for computing and plotting CDFs for Models I & III when parameter $\gamma = \frac{m+1}{2}$:

| | |
|---|---|
| $\operatorname{erf}(x) \;= \Phi(x) = \frac{2}{\sqrt{\pi}} \int_0^x e^{-t^2}\, dt = \frac{\gamma\left(\frac{1}{2}, x^2\right)}{\sqrt{\pi}}$ | Definition in GR 8.250.1 |
| $\operatorname{erf}(-\infty) = \Phi(-\infty) = -1$ <br> ▲may be questionable but follows from the math of the (even-function) integrand of the integral | $\operatorname{erf}(-\infty) \;= \frac{2}{\sqrt{\pi}} \int_0^{-\infty} e^{-t^2}\, dt = -\frac{2}{\sqrt{\pi}} \int_{-\infty}^{0} e^{-t^2}\, dt = -\frac{2}{\sqrt{\pi}} \int_0^{\infty} e^{-t^2}\, dt = -1$ <br> (even-function property on $\pm\infty$ ▼) <br> 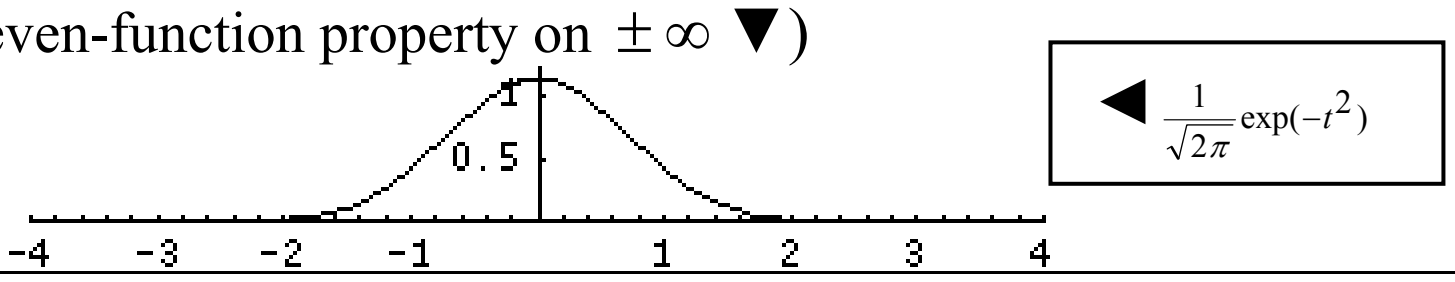  |
| $\operatorname{erf}(0) \;= \Phi(0) = 0$ | $\frac{2}{\sqrt{\pi}} \int_0^0 e^{-t^2}\, dt = 0$ |
| $\operatorname{erf}(+\infty) = \Phi(+\infty) = 1$ | $\frac{2}{\sqrt{\pi}} \int_0^{\infty} e^{-t^2}\, dt = 1$ |

- Numerical methods are required to compute the Model I CDF $F(x)$ for $\gamma \neq \frac{m+1}{2}$ (vs. $\gamma = \frac{m+1}{2}$ when the error function erf is easily implemented as in above CDF plots). One approach to compute $\frac{1}{2}\left[1 + \frac{\gamma(\nu, \beta x^n)}{\Gamma(\nu)}\right]$ may be done by a truncated series of an infinite expansion, GR 8.354.1,

$$\gamma(a, x) = \frac{\sum_{n=0}^{\infty} (-1)^n x^{a+n}}{n!(a+n)}$$

For example, to compute, $\frac{1}{2}\left[1 + \frac{\gamma\left(\frac{1}{4}, \frac{x^4}{2}\right)}{\Gamma\left(\frac{1}{4}\right)}\right]$,

$$F(x) = \frac{1}{2}\left[1 + \frac{\gamma\left(\frac{1}{4}, \frac{x^4}{2}\right)}{\Gamma\left(\frac{1}{4}\right)}\right] = \frac{1}{2}\left[1 + \frac{1}{\Gamma\left(\frac{1}{4}\right)} \frac{4\sum_{n=0}^{\infty} (-1)^n x^{4n+1}}{n!\left(2^{n+\frac{1}{4}}\right)(4n+1)}\right]$$

for the range of $x$ of say, $-6 < x < 6$, to a predetermined number of iterations for convergence (**TBD**)

Figure A2a.

# Complementary CDF Plots

<table>
<tr><th>MODEL I<br>(Even Function, Symmetric Density)</th><th>MODEL II</th></tr>
<tr><td>

**1-*F*(*x*)**

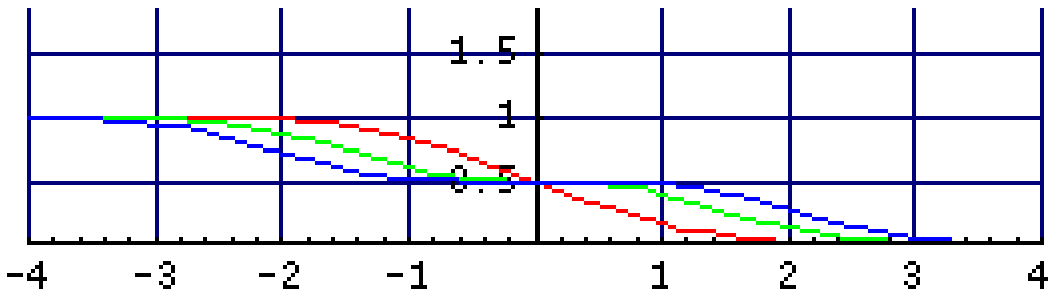

</td><td>

**1-*F*(*x*)**

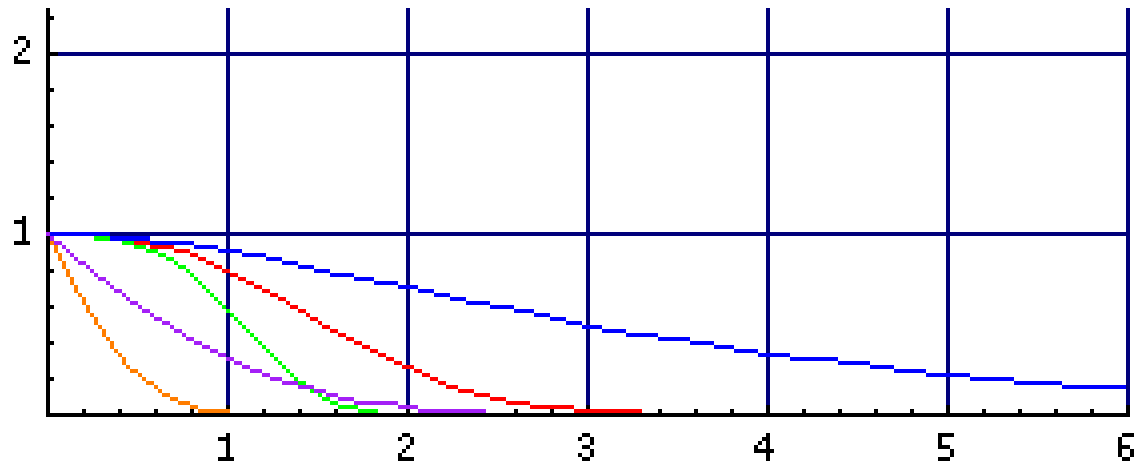

</td></tr>
<tr><td>▲Plot Parameters  Same as CDF</td><td>▲Plot Parameters  Same as CDF</td></tr>
</table>

NOTEs:

- $0 \le 1 - F(x) \le 1$.
- The complementary CDF, $1 - F(x)$, is also called the Reliability or Survival Function
- GEMI, $F(-x) = 1 - F(x)$

Figure A2b.

# Hazard Function (Failure Rate) Plots

| **MODEL I**<br>**(Even Function, Symmetric Density)** | **MODEL II** |
|---|---|
| **N/A**<br><br>**Since failure rate function** $h(t)$ **is defined on time,** $0 < t < \infty$**,**<br>**it is not applicable for**<br>$\pm\infty$ **Models I and III**<br>(Hoel, et al., Vol. I, p. 137, EX. 40) | ***h(t)***<br>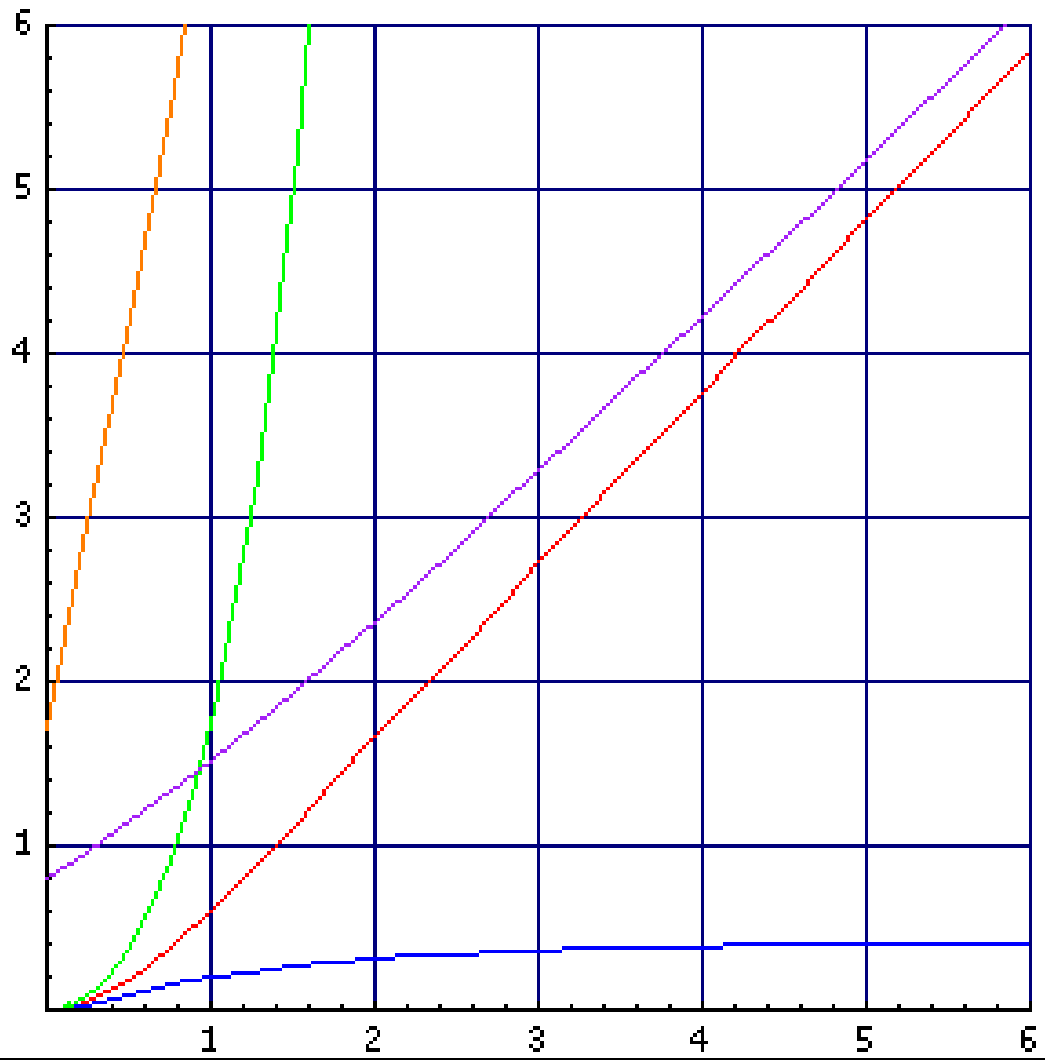<br>▲Plot Parameters same as CDF<br><br>• Blue is nearly a constant low failure-rate function.<br>• The Red (Maxwell) and Purple (doubled half-Gaussian) functions are linear<br>• Green and Yellow are very steep with fast, high failure rates |

**NOTE:** Hazard function (failure rate), $h(t)$, is the tendency for a system to fail in the next time interval $t$ given that the system has not failed prior to time $t$.

$$h(t) = \lim_{\Delta t \to 0} \frac{\Pr\left(t \le T \le t + \Delta t \middle| T > t\right)}{\Delta t} = \frac{f(t)}{1 - F(t)}, 0 < t < \infty,\ f(0) \le h(t) < \infty$$

- EXAMPLES
  - GEMII (all models), $h(t) = \dfrac{n\beta^{\gamma} t^{m}}{\Gamma\left(\gamma, \beta t^{n}\right)} e^{-\beta t^{n}}$ (variable rate)
  - GEMII, Exponential, $h(t) = \lambda$ (constant rate)
  - GEMII, Weibull, $h(t) = abt^{b-1}$ (variable rate)

# MODEL II PDF DISTRIBUTIONS

**Table A3. Model Comparison for Exponential-Gamma Distributions**

We can examine the structural similarities among the interrelated distributions—Exponential, Gamma, Weibull, transformed gamma distribution, generalized gamma, and Model II—by comparing their PDFs using a standard notation ~ Klugman et al., 2004 (Appendix A).

- Setting $k=1$ reveals an interesting set of distributions of the exponential family

| Distribution Model (Parameters) | PDF, $f(x)$ | PDF, $f(x)$ $k=1$ | Shape Parameter | Scale Parameter |
|---|---|---|---|---|
| **Exponential** $(\theta)$ | $\frac{e^{-\frac{x}{\theta}}}{\theta}$ | $\frac{e^{-\frac{x}{\theta}}}{\theta}$ | — | $\theta$ |
| **Gamma** $(k,\theta)$ | $\frac{\left(\frac{x}{\theta}\right)^k}{x\Gamma(k)}e^{-\frac{x}{\theta}}$ | $\frac{e^{-\frac{x}{\theta}}}{\theta}$ | $k$ | $\theta$ |
| **Weibull** $(k,\theta)$ | $\frac{k\left(\frac{x}{\theta}\right)^k}{x}e^{-\left(\frac{x}{\theta}\right)^k}$ | $\frac{e^{-\frac{x}{\theta}}}{\theta}$ | $k$ | $\theta$ |
| **Generalized Gamma** (Stacy, 1962) $(k,\theta,d)$ | $k\frac{\left(\frac{x}{\theta}\right)^d}{x\Gamma\left(\frac{d}{k}\right)}e^{-\left(\frac{x}{\theta}\right)^k}$ | $\frac{x^{d-1}}{\theta^d\Gamma(d)}e^{-\frac{x}{\theta}}$ | $k,d$ | $\theta$ |
| **Transformed Gamma** (Venter, 1983) $(k,\theta,r)$ | $k\frac{\left(\frac{x}{\theta}\right)^{kr}}{x\Gamma(r)}e^{-\left(\frac{x}{\theta}\right)^k}$ | $\frac{x^{r-1}}{\theta^r\Gamma(r)}e^{-\frac{x}{\theta}}$ | $k,r$ | $\theta$ |
| **GEM II** $(k,\theta,\gamma)$ $\gamma=\frac{m+1}{n}=\frac{m+1}{k}$, $m=n\gamma-1=k\gamma-1$ | $k\frac{\left(\frac{x}{\theta^{\frac{1}{k}}}\right)^{k\gamma}}{x\Gamma(\gamma)}e^{-\left[\theta^{k-1}\left(\frac{x}{\theta}\right)^k\right]}$ | $\frac{x^{\gamma-1}}{\theta^\gamma\Gamma(\gamma)}e^{-\frac{x}{\theta}}$ $\gamma=m+1$ | $k,\gamma$ | $\theta$ |

- NOTE: Scale parameters are called rate or inverse scale parameters in cases where they appear in numerators although some authors convert scale parameters to rate parameters by re-writing the PDFs, as in the exponential distribution, sometimes written in a form such as,
  - $f(x)=\lambda e^{-\lambda x}$, where $\lambda=\frac{1}{\theta}$

- If $k=\gamma=1$ in Model II, the PDF $f(x;k=\gamma=1)=\frac{e^{-\frac{x}{\theta}}}{\theta}$, the same as the exponential family of ~ exponential, gamma, Weibull. If $r=d=1$, the Transformed Gamma & Generalized Gamma also reduce to the same. Hence, the suggested model-name—Generalized Exponential Distribution I,…, Generalized Exponential Distribution IV.

- Stacy's and Venter's 3-parameter models are generalizations of the gamma distribution derived from the following modeling functions, $g(x)$, and are related functionally by $d = kr$.
    - For Stacy's generalized gamma distribution, $g(x) = x^{d-1}\exp\left(-\frac{x}{\theta}\right)^{k}$
    - For Venter's transformed gamma distribution, $g(x) = x^{kr-1}\exp\left(-\frac{x}{\theta}\right)^{k}$
        - GEMII—the generalized exponential distribution—differs since it is a generalization of all similar elementary algebraic-exponential probability functions on $[0,\infty)$ derived from the basis or modeling function, $g(x) = x^{k\gamma-1}\exp\left[-\theta^{k-1}\left(\frac{x}{\theta}\right)^{k}\right], \gamma = \frac{m+1}{k}$
            - Only when $k = 1$ are all distributions of the same structure, with shape parameter, $d = r = m+1$.

# Change of Variable Distributions

**Table A4. Summary of Probability Models for Selected Transforms**

## PDFs

| MODEL | DOMAIN PDF, $f(x)$ | TRANSFORM | PDF $g(y)$ | PARAMETER RESTRICTION |
|---|---|---|---|---|
| I | $-\infty < x < \infty$ | $Y = a + bX$ | $\frac{n\beta^{\gamma}\left(\frac{y-a}{b}\right)^{m}}{2\lvert b\rvert\Gamma(\gamma)} e^{-\beta\left(\frac{y-a}{b}\right)^{n}}, -\infty < y < \infty$ | $z > 0$<br>$\gamma = \frac{m+1}{n} > 0$<br>$m > -1, \beta, n > 0$<br>$c, k \neq 0$<br>symmetric density |
| | | $z = \lvert x\rvert$<br>$Y = Z^{c}$ | $\frac{n\beta^{\gamma} y^{\frac{n\gamma}{c}-1}}{\lvert c\rvert\Gamma(\gamma)} e^{-\beta y^{\frac{n}{c}}}, y > 0$ | |
| | | $z = \lvert x\rvert$<br>$Y = (kZ)^{c}$ | $\frac{n\beta^{\gamma} y^{\frac{n\gamma}{c}-1}}{\lvert ck\rvert\Gamma(\gamma)k^{m}} e^{-\frac{\beta y^{\frac{n}{c}}}{k^{n}}}, y > 0$ | |
| | | $z = \lvert x\rvert$<br>$Y = \frac{1}{(kZ)^{c}}$ | $\frac{n\beta^{\gamma}}{\lvert ck\rvert\Gamma(\gamma)k^{m} y^{\frac{n\gamma}{c}+1}} e^{-\frac{\beta}{k^{n}}\frac{1}{y^{\frac{n}{c}}}}, y > 0$ | |
| II | $0 \leq x < \infty$ | $Y = X^{c}$ | $\frac{n\beta^{\gamma} y^{\frac{n\gamma}{c}-1}}{\lvert c\rvert\Gamma(\gamma)} e^{-\beta y^{\frac{n}{c}}}, y > 0$ | $\gamma = \frac{m+1}{n} > 0$<br>$m > -1, \beta, n > 0$<br>$c, k \neq 0$ |
| | | $Y = \frac{1}{X}$ | $\frac{n\beta^{\gamma}}{\Gamma(\gamma) y^{n\gamma+1}} e^{-\frac{\beta}{y^{n}}}, y > 0$ | |
| | | $Y = (kX)^{c}$ | $\frac{n\beta^{\gamma} y^{\frac{n\gamma}{c}-1}}{\lvert ck\rvert\Gamma(\gamma)k^{m}} e^{-\frac{\beta y^{\frac{n}{c}}}{k^{n}}}, y > 0$ | |
| | | $Y = \frac{1}{(kX)^{c}}$ | $\frac{n\beta^{\gamma}}{\lvert ck\rvert\Gamma(\gamma)k^{m} y^{\frac{n\gamma}{c}+1}} e^{-\frac{\beta}{k^{n}}\frac{1}{y^{\frac{n}{c}}}}, y > 0$ | |
| | | $Y = a + bX$ | $\frac{n\beta^{\gamma}\left(\frac{y-a}{b}\right)^{m}}{\lvert b\rvert\Gamma(\gamma)} e^{-\beta\left(\frac{y-a}{b}\right)^{n}}, y > 0$ | |
| III | $-\infty < x < \infty$ | $Y = a + bX$ | | |
| | | $z=\left\lvert\frac{x-a}{b}\right\rvert$ <br> $Y=Z^c$ | | $z>0$ <br> $\gamma=\frac{m+1}{n}>0$ <br> $m>-1,\beta,b,n>0$ <br> $-\infty<a<\infty$ <br> $c,k\neq 0$ <br> symmetric density |
| | | $z=\left\lvert\frac{x-a}{b}\right\rvert$ <br> $Y=(kZ)^c$ | Same as I | |
| | | $z=\left\lvert\frac{x-a}{b}\right\rvert$ <br> $Y=\frac{1}{(kZ)^c}$ | | |
| IV | $a\leq x<\infty$ | $Y=\left(\frac{X-a}{b}\right)$ | Same as II | $\gamma=\frac{m+1}{n}>0$ <br> $m>-1,\beta,b,n>0$ <br> $-\infty<a<\infty$ <br> $c,k\neq 0$ |
| | | $Y=\frac{1}{\left(\frac{X-a}{b}\right)^c}$ | | |
| | | $Y=k\left(\frac{X-a}{b}\right)$ | | |
| | | $Y=\frac{1}{k\left(\frac{X-a}{b}\right)^c}$ | | |

NOTE: Transformed by Change of Variable, Theorem 1, Hoel et al., Vol. I, p.119, $g(y)=f(x)\left\lvert\frac{dx}{dy}\right\rvert$.

For GEMI symmetric PDF on $\pm\infty$, first set $z=\lvert x\rvert, z>0$ in $f(x)$. Find new PDF $f(z)$ on domain $0\to\infty$, $f(z)=\frac{n\beta^\gamma z^m}{\Gamma(\gamma)}e^{-\beta z^n}$, and transform by $Y=Z^c$ to obtain transformed density,

$$g(y)=f(z)\left\lvert\frac{dz}{dy}\right\rvert=f(y^{\frac{1}{c}})\left\lvert\frac{dy^{\frac{1}{c}}}{dy}\right\rvert=\frac{n\beta^\gamma y^{\frac{n\gamma}{c}-1}}{\lvert c\rvert\Gamma(\gamma)}e^{-\beta y^{\frac{n}{c}}},\ y>0.$$

Check $\int_0^\infty g(y)dy=1$ by GR 3.326.2. Other GEMI transforms are similar. GEMIII transforms follow same logic for $z=\left\lvert\frac{x-a}{b}\right\rvert>0$.

- Special treatment for densities on $\pm\infty$ is required by Theorem I for transform $Y$ to be strictly increasing or decreasing and differentiable on domain of density. Cf. J.E. Freund, *Math. Stat.*, 2nd ed., pp. 123-4 & Bers, pp. 93-4 ("Monotone Functions").
- GEMII & IV PDFs follow from Theorem 1 by selected transforms.

- Example (Hoel, et al., Vol I, p. 136, EX. 32); find density $Y=|X|$ for $n\left(0,\sigma^2\right)$.

  $f(x)=\dfrac{e^{-\frac{1}{2}\left(\frac{x}{\sigma}\right)^2}}{\sigma\sqrt{2\pi}}$ ;for GEMIII solution with parameters

  $\left[m=0,n=2,\beta=\dfrac{1}{2\sigma^2},c=k=1\right], g(y)=\dfrac{2}{|1|\Gamma\left(\frac{1}{2}\right)}\left(\dfrac{1}{2\sigma^2}\right)^{\frac{1}{2}}e^{-\frac{1}{2}\left(\frac{y}{\sigma}\right)^2}=\dfrac{2}{\sigma\sqrt{2\pi}}e^{-\frac{1}{2}\left(\frac{y}{\sigma}\right)^2}, y>0$

  - Check calculation by GR 3.326.2 to show, $\int_0^\infty g(y)dy=1$.

- Example; GEMI PDF for $Y=Z^2$ for $m=0,n=2,\beta=\gamma=1/2$, $g(y)=\dfrac{y^{-\frac{1}{2}}}{\sqrt{2\pi}}e^{-\frac{y}{2}}$ (Gaussian)

- Example; GEMI PDF, $f(x)=\dfrac{x^4}{3\sqrt{2\pi}}e^{-\frac{x^2}{2}}$. Find density $Y=\dfrac{1}{Z^2}; g(y)=\dfrac{e^{-\frac{1}{2y}}}{3\sqrt{2\pi}y^{\frac{7}{2}}}$ by

  $Y=\dfrac{1}{(kZ)^2}$ for parameters $[m=4,n=2,\beta=1/2,\gamma=5/2,k=1,c=2]$. Check by integral in NOTE ▼

- Example; for PDF $f(x)=\dfrac{2x^2}{\sqrt{2\pi}}e^{-\frac{x^2}{2}}$ find GEMII PDF $Y=X^{-1}$;

  $g(y)=\dfrac{2e^{-\frac{1}{2y^2}}}{y^4\sqrt{2\pi}}\left(m=n=2;\beta=\dfrac{1}{2};\gamma=\dfrac{3}{2}\right)$. Check by integral in NOTE ▼

- Example (Hoel, et al., Vol. I, pp. 120-1);Find density $Y=X^{\frac{1}{\beta}}$ for exponential density $\lambda e^{-\lambda x}$;
  $g(y)=|\beta|\lambda y^{\beta-1}e^{-\lambda y^\beta}$

- Example; find GEMII density for $Y=\dfrac{1}{\beta X^n}$; use transform $Y=\dfrac{1}{(kX)^c}=\left(\beta^{\frac{1}{n}}X\right)^{-n}$ to get,

  $g(y)=\dfrac{e^{-\frac{1}{y}}}{y^{\gamma+1}\Gamma(\gamma)}$. Check by integral in NOTE ▼

- Example; GEMII gamma PDF for $Y=\dfrac{1}{X}$; $g(y)=\dfrac{\lambda^p}{\Gamma(p)y^{p+1}}e^{-\frac{\lambda}{y}}$
  - Exponential; Put $p=1$ for $g(y)$

- Example; GEMII Chi-square PDF for $Y=\dfrac{1}{X}$; $g(y)=\dfrac{e^{-\frac{1}{2y}}}{2^{\frac{\nu}{2}}\Gamma\left(\frac{\nu}{2}\right)y^{\frac{\nu}{2}+1}}$ with parameters

  $\left[m=\dfrac{\nu}{2}-1,n=1,\beta=\dfrac{1}{2},\gamma=\dfrac{\nu}{2}\right]$

- NOTE: the following integral from Table A7 is useful for inverse change of variable transformed distributions:

$$\int_u^v \frac{e^{-\frac{\beta}{x^n}}}{x^m}dx = \frac{\Gamma\left(\frac{m-1}{n},\frac{\beta}{v^n}\right)-\Gamma\left(\frac{m-1}{n},\frac{\beta}{u^n}\right)}{n\beta^{\frac{m-1}{n}}}$$ by GR 2.33.19 with $n \to -n$

  - Example: the GEMII PDF $g(y)$ for inverse $Y = \frac{1}{X}$, can be shown to equal 1 by setting

$m \to n\gamma + 1$, and $\begin{cases} \left\{ \lim\limits_{v\to\infty} \frac{\beta}{v^n} = 0 \right. \\ \left\{ \lim\limits_{u\to 0} \frac{\beta}{u^n} = \infty \right. \end{cases}$ and, $\begin{cases} \Gamma(a,0) = \Gamma(a) \\ \Gamma(a,\infty) = 0 \end{cases}$ , to give integral *I*,

$$I = \frac{n\beta^\gamma}{\Gamma(\gamma)}\int_0^\infty \frac{e^{-\frac{\beta}{y^n}}}{y^{n\gamma+1}}dy = \frac{n\beta^\gamma}{\Gamma(\gamma)}\left(\frac{\Gamma(\gamma,0)-\Gamma(\gamma,\infty)}{n\beta^\gamma}\right) = 1$$

- Thus, $\int_0^\infty \frac{e^{-\frac{\beta}{x^n}}}{x^m}dx = \frac{\Gamma\left(\frac{m-1}{n}\right)}{n\beta^{\frac{m-1}{n}}}, \frac{m-1}{n} > 0$

$$\int_0^u \frac{e^{-\frac{\beta}{x^n}}}{x^m}dx = \frac{\Gamma\left(\frac{m-1}{n},\frac{\beta}{u^n}\right)}{n\beta^{\frac{m-1}{n}}}$$

$$\int_u^\infty \frac{e^{-\frac{\beta}{x^n}}}{x^m}dx = \frac{\Gamma\left(\frac{m-1}{n}\right)-\Gamma\left(\frac{m-1}{n},\frac{\beta}{u^n}\right)}{n\beta^{\frac{m-1}{n}}} = \frac{\gamma\left(\frac{m-1}{n},\frac{\beta}{u^n}\right)}{n\beta^{\frac{m-1}{n}}}, \frac{m-1}{n} > 0$$

# Moments

$$EY^j = \int_0^\infty y^j g(y)dy$$

| MODEL | TRANSFORM | $EY^j$ | PARAMETER RESTRICTION |
|---|---|---|---|
| I | $Y = Z^c$ | $\dfrac{\Gamma\left(\gamma + j\frac{c}{n}\right)}{\Gamma(\gamma)\beta^{j\frac{c}{n}}}$ | $z = \lvert x \rvert, z > 0$ <br> $\gamma = \dfrac{m+1}{n} > 0$ <br> $m > -1, \beta, n > 0$ <br> $c, k \neq 0$ <br> $\gamma \pm j\dfrac{c}{n} > 0$ |
| | $Y = (kZ)^c$ | $\left(\dfrac{k}{\beta^{\frac{1}{n}}}\right)^{cj} \dfrac{\Gamma\left(\gamma + j\frac{c}{n}\right)}{\Gamma(\gamma)}$ | |
| | $Y = \dfrac{1}{(kZ)^c}$ | $\left(\dfrac{\beta^{\frac{1}{n}}}{k}\right)^{cj} \dfrac{\Gamma\left(\gamma - j\frac{c}{n}\right)}{\Gamma(\gamma)}$ | |
| II | $Y = X^c$ | $\dfrac{\Gamma\left(\gamma + j\frac{c}{n}\right)}{\Gamma(\gamma)\beta^{j\frac{c}{n}}}$ | $\gamma = \dfrac{m+1}{n} > 0$ <br> $m > -1, \beta, n > 0$ <br> $c, k \neq 0$ <br> $\gamma \pm j\dfrac{c}{n} > 0$ |
| | $Y = \dfrac{1}{X^c}$ | $\dfrac{\beta^{j\frac{c}{n}}\Gamma\left(\gamma - j\frac{c}{n}\right)}{\Gamma(\gamma)}$ | |
| | $Y = (kX)^c$ | $\left(\dfrac{k}{\beta^{\frac{1}{n}}}\right)^{cj} \dfrac{\Gamma\left(\gamma + j\frac{c}{n}\right)}{\Gamma(\gamma)}$ | |
| | $Y = \dfrac{1}{(kX)^c}$ | $\left(\dfrac{\beta^{\frac{1}{n}}}{k}\right)^{cj} \dfrac{\Gamma\left(\gamma - j\frac{c}{n}\right)}{\Gamma(\gamma)}$ | |
| III | $Y = (kZ)^c$ | Same as I | $z = \left\lvert \dfrac{x-a}{b} \right\rvert, z > 0$ <br> $\gamma = \dfrac{m+1}{n} > 0$ <br> $m > -1, \beta, b, n > 0$ <br> $-\infty < a < \infty$ <br> $c, k \neq 0$ <br> $\gamma \pm j\dfrac{c}{n} > 0$ |
| | $Y = (kZ)^c$ | | |
| | $Y = \dfrac{1}{(kZ)^c}$ | | |
| IV | $Y=\left(\frac{X-a}{b}\right)^c$ | Same as II | $\gamma=\frac{m+1}{n}>0$<br>$m>-1,\beta,b,n>0$<br>$-\infty<a<\infty$<br>$c,k\neq 0$<br>$\gamma\pm j\frac{c}{n}>0$ |
| | $Y=\frac{1}{\left(\frac{X-a}{b}\right)^c}$ | | |
| | $Y=\frac{1}{k\left(\frac{X-a}{b}\right)}$ | | |

- NOTE: integrals derived from either GR 3.326.2, or integral in NOTE under PDFs, Table A4.
- Example, GEMI mean for $Y=Z^2$ if $m=0,n=2,$ $\mu=\frac{\gamma}{\beta};\sigma^2=\frac{\gamma}{\beta^2}$;e.g., normal transformed dist. mean is $1,\sigma^2=2$.
- Example: GEMII transformed gamma mean for $Y=\frac{1}{X}$, $\mu=\frac{\lambda\Gamma(p-1)}{\Gamma(p)}=\frac{\lambda}{p-1}$ with parameters $[n=1,m=p-1,\gamma=p,\beta=\lambda]$
  - Since exponential dist. is gamma with $p=1$, exp. moments do not exist with this transform; e.g.,

$$\mu=\lambda\int_0^\infty\frac{e^{-\frac{\lambda}{y}}}{y}dy=\lambda\Gamma(0)\to\infty$$

- Example; GEMIII mean for density $Y=\left|\frac{x-a}{b}\right|^2$ if $f(x)=\frac{1}{10\sqrt{2\pi}}\left(\frac{x-5}{10}\right)^2\exp\frac{1}{2}\left(\frac{x-5}{10}\right)^2$;

  $\mu=1$

# CDFs

$$F(y)=\int_0^y g(t)dt$$

| MODEL | TRANSFORM | $F(y)$ | PARAMETER RESTRICTION |
|---|---|---|---|
| I | $Y=Z^c$ | $\frac{\gamma\left(\nu,\beta y^{\frac{n}{c}}\right)}{\Gamma(\gamma)}$ | $z=\lvert x\rvert, z>0$<br>$\gamma=\nu=\frac{m+1}{n}$<br>$m>-1,\beta,n>0$<br>$c,k\neq 0$ |
| | $Y=(kZ)^c$ | $\frac{\gamma\left(\nu,\frac{\beta}{k^n}y^{\frac{n}{c}}\right)}{\Gamma(\gamma)}$ | |
| | $Y=\frac{1}{(kZ)^c}$ | $\frac{\Gamma\left(\gamma,\frac{\beta/k^n}{y^{\frac{n}{c}}}\right)}{\Gamma(\gamma)}$ | |
| II | $Y=X^c$<br>&<br>$Y=\frac{1}{X^c}$ | $\begin{cases}\frac{\gamma\left(\nu,\beta y^{\frac{n}{c}}\right)}{\Gamma(\gamma)}, c>0\\ \frac{\Gamma\left(\gamma,\frac{\beta}{y^{\frac{n}{c}}}\right)}{\Gamma(\gamma)}, c<0\end{cases}$ | $\gamma=\nu=\frac{m+1}{n}$<br>$m>-1,\beta,n>0$<br>$c,k\neq 0$ |
| | $Y=(kX)^c$ | $\frac{\gamma\left(\nu,\frac{\beta}{k^n}y^{\frac{n}{c}}\right)}{\Gamma(\gamma)}$ | |
| | $Y=\frac{1}{(kX)^c}$ | $\frac{\Gamma\left(\nu,\frac{\beta/k^n}{y^{\frac{n}{c}}}\right)}{\Gamma(\gamma)}$ | |
| III | $Y=Z^c$ | Same as I | $z=\left\lvert\frac{x-a}{b}\right\rvert>0$<br>$\gamma=\nu=\frac{m+1}{n}$<br>$m>-1,\beta,b,n>0$<br>$-\infty<a<\infty$<br>$c,k\neq 0$ |
| | $Y=(kZ)^c$ | | |
| | $Y=\frac{1}{(kZ)^c}$ | | |
| IV | $Y=\left(\frac{X-a}{b}\right)^{c}$ | Same as II | $\gamma=\nu=\frac{m+1}{n}$<br>$m>-1,\beta,b,n>0$<br>$-\infty<a<\infty$<br>$c,k\neq 0$ |
| | $Y=1\Big/\left(\frac{X-a}{b}\right)^{c}$ | | |

- NOTE: integrals derived from either GR 3.326.2, or integral in NOTE under PDFs, Table A4.
- NOTE alternative expressions possible by GR 8.356.3, $\Gamma(\alpha)=\gamma(\alpha,x)+\Gamma(\alpha,x)$

- Example GEMI normal for $Y=Z^2$, $F(y)=\dfrac{\gamma\left(\frac{1}{2},\frac{y}{2}\right)}{\Gamma\left(\frac{1}{2}\right)}=\mathrm{erf}\left(\sqrt{\frac{y}{2}}\right)$.

- Example: GEMII, Weibull for $Y=\dfrac{1}{X}$, $F(y)=\dfrac{\Gamma\left(1,\frac{a}{y^b}\right)}{\Gamma(1)}=e^{-\frac{a}{y^b}}$ with parameters $[n=b,m=b-1,\gamma=1,\beta=a]$

# Log-Model Distributions

**Table A5. Summary of Log Probability Models**

## Log PDFS

| MODEL | DOMAIN PDF, $f(x)$ | LOG PDF $g(y)$ | PARAMETER RESTRICTION |
|---|---|---|---|
| I | $-\infty < x < \infty$ | $\frac{n\beta^{\gamma}(\ln y)^{m} e^{-\beta(\ln y)^{n}}}{2y\Gamma(\gamma)}$, $0 < y < \infty$ | $Y = e^{X}$<br>$\gamma = \frac{m+1}{n} > 0$<br>$m > -1, \beta, n > 0$<br>Symmetric density |
| II | $0 \le x < \infty$ | $\frac{n\beta^{\gamma}(\ln y)^{m} e^{-\beta(\ln y)^{n}}}{y\Gamma(\gamma)}$, $1 < y < \infty$ | $Y = e^{X}$<br>$\gamma = \frac{m+1}{n} > 0$<br>$m > -1, \beta, n > 0$ |
| III | $-\infty < x < \infty$ | $\frac{n\beta^{\gamma}\left(\frac{\ln y - a}{b}\right)^{m} e^{-\beta\left(\frac{\ln y - a}{b}\right)^{n}}}{2by\Gamma(\gamma)}$, $0 < y < \infty$ | $Y = e^{X}$<br>$\gamma = \frac{m+1}{n} > 0$<br>$m > -1, \beta, n, b > 0$<br>$-\infty < a < \infty$<br>Symmetric density |
| IV | $a \le x < \infty$ | $\frac{n\beta^{\gamma}\left(\frac{\ln y - a}{b}\right)^{m} e^{-\beta\left(\frac{\ln y - a}{b}\right)^{n}}}{by\Gamma(\gamma)}$, $e^{a} < y < \infty$ | $Y = e^{X}$<br>$\gamma = \frac{m+1}{n} > 0$<br>$m > -1, \beta, n, b > 0$<br>$-\infty < a < \infty$ |

- Example, GEM I, Change of Variable (Theorem 1, Hoel, et al., Vol. I, p. 119)

$$Y = e^{X}, X = \ln Y, \left|\frac{dx}{dy}\right| = \left|\frac{d\ln y}{dy}\right| = \left|\frac{1}{y}\right|; Y_{-\infty} = e^{-\infty} = 0; Y_{\infty} = e^{\infty} = \infty;\ g(y) = f(x)\left|\frac{dx}{dy}\right| = f(\ln y)\left|\frac{d\ln y}{dy}\right|$$

  - GEM II $g(y)$ on $x \ge 0$ is similar with limits $Y_0 = 1, Y_{\infty} = \infty$.
    - GEM III & IV similar with limits $(Y_{-\infty}, Y_{\infty})$ & $(Y_a, Y_{\infty})$, respectively.

- Example, lognormal, GEM I: $g(y) = \frac{1}{y\sqrt{2\pi}} e^{-\frac{1}{2}(\ln y)^2}$ with parameters $[m = 0, n = 2, \beta = \gamma = 1/2]$

- Example, loggamma, GEM II: $g(y) = \frac{\lambda^{p}(\ln y)^{p-1} e^{-\lambda \ln y}}{y\Gamma(p)} = \frac{\lambda^{p}(\ln y)^{p-1}}{\Gamma(p) y^{\lambda+1}}$ with parameters $[m = p-1, n = 1, \beta = \lambda, \gamma = p]$

- Example, GEM III lognormal: $g(y) = \dfrac{e^{-\frac{1}{2}\left(\frac{\ln y-\mu}{\sigma}\right)^2}}{y\sigma\sqrt{2\pi}}$ with parameters $\left[m = 0, n = 2, \beta = \gamma = 1/2, a = \mu, b = \sigma\right]$

- Example, the mode of the lognormal, GEM III. From the definition given in Notes of Table 3c, the mode is defined as the derivative of the PDF evaluated at 0; i.e., $\left(\dfrac{d}{dy}g(y)\right)_{y=0}$ for the PDF, $g(y) = \dfrac{1}{y\sigma\sqrt{2\pi}}e^{-\frac{1}{2}\left(\frac{\ln y-\mu}{\sigma}\right)^2}$, so that taking the derivative of the PDF, setting it equal to 0 and solving for $y$ gives: $\left(\dfrac{d}{dy}g(y)\right)_{y=0} = \dfrac{d}{dy}\dfrac{1}{y\sigma\sqrt{2\pi}}e^{-\frac{1}{2}\left(\frac{\ln y-\mu}{\sigma}\right)^2} = 0$ which after some manipulation provides the mode Mo, $\text{Mo} = y = e^{\mu-\sigma^2}$

- Example, GEMIV loggamma: $g(y) = \dfrac{\lambda^p(\ln y-\gamma)^{p-1}}{y\Gamma(p)}e^{-\lambda(\ln y-\gamma)}$ with parameters $\left[m = p-1, n = 1, \beta = \lambda, a = \gamma, b = 1\right]$

NOTE: A reverse log transform can be made by the transforms (Lawless, 1980): $X = e^Y; Y = \ln X; \beta = e^U; n = \frac{1}{n}; \left|\frac{dx}{dy}\right| = e^Y; Y_0 = -\infty, Y_\infty = +\infty$ . Then for Model II, we obtain the transformed PDF using rate parameter, $\beta \to \beta^{-1}$, i.e., $f(x) = \dfrac{n}{\beta^\gamma\Gamma(\gamma)}x^{n\gamma-1}\exp\left(-\dfrac{x^n}{\beta}\right)$, where $m = n\gamma - 1$,

$$g(y) = f\left(e^Y\right)\left|\frac{de^y}{dy}\right| = \frac{1}{n\Gamma(\gamma)}\frac{1}{\left(e^u\right)^\gamma}\left(e^y\right)^{\frac{\gamma}{n}-1}\exp\left[\frac{\left(e^y\right)^{\frac{1}{n}}}{e^u}\right]e^y = \frac{1}{n\Gamma(\gamma)}\exp\left[\gamma\frac{y-un}{n} - \exp\left(\frac{y-un}{n}\right)\right], -\infty < y < \infty$$

This model is useful for maximum likelihood estimation (Table A6). See Lawless, 1980 or Bell, 1988 for application to the generalized gamma distribution (Stacy, 1962).

NOTE: still other transforms for Model II can be obtained from the relation, $\beta X^n = -\phi\ln Y$, where $\phi$ is a "fake parameter". Isolate the quantities $X, Y, \left|\frac{dx}{dy}\right|$, which results in the PDF called the log-gamma distribution, which has been used for likelihood ratio test (LRT) of the Model II Stacy generalized gamma distribution, $g(y) = \dfrac{\varphi^\gamma}{\Gamma(\gamma)}y^{\varphi-1}(-\ln y)^{\gamma-1}, 0 < y < 1, \text{Re}\,\phi, \gamma > 0$ **.** See Consul & Jain, 1971. Verify integral by a change of variable.

NOTE: Using Crooks limit transformation on the Model IV Amoroso distribution, the Model IV PDF can be written as an exp. transform**,**

$f(x) = \dfrac{\beta^{\frac{1}{n}}\gamma^\gamma}{b\Gamma(\gamma)}\exp\left\{\gamma\left(\beta^{\frac{1}{n}}\dfrac{x-a}{b}\right) - \gamma\exp\left[\beta^{\frac{1}{n}}\dfrac{x-a}{b}\right]\right\}, \gamma = \dfrac{m+1}{n}, -\infty < x < \infty$ **.** For Model II, **set** $a = 0, b = 1$. This formulation is good for many distributions, including non-Gaussian and skewed distributions.

# Log Moments

$$EY^j = \int_{-\infty}^{\infty} y^j g(y) dy$$

| MODEL | MOMENTS FUNCTION $EY^j$ | NOTES |
|---|---|---|
| I | $\frac{n\beta^\gamma}{2\Gamma(\gamma)} \int_{-\infty}^{\infty} s^m e^{-(\beta s^n - js)} ds$ <br> $s = \ln y$ | $\gamma = \frac{m+1}{n} > 0,\ m > -1, \beta, n > 0$ <br> EX: lognormal by GR 3.323.2, for $n = 2, m = 0, \beta = \frac{1}{2}$, <br> $\int_{-\infty}^{\infty} \exp(-p^2x^2 \pm qx)dx = \exp\left(\frac{q^2}{4p^2}\right)\frac{\sqrt{\pi}}{p}; p^2 = \frac{1}{2}; q = j$ <br> $EY^j = e^{\frac{j^2}{2}}$ |
| II | $\frac{n\beta^\gamma}{\Gamma(\gamma)} \int_{0}^{\infty} s^m e^{-(\beta s^n - js)} ds$ <br> $s = \ln y$ | $\gamma = \frac{m+1}{n} > 0,\ m > -1, \beta, n > 0$ <br> EX: loggamma mean for $n = 1, m = p - 1, \beta = \lambda$, <br> $\mu = \left(\frac{\lambda}{\lambda - 1}\right)^p$ from $EY^j = \left(\frac{\lambda}{\lambda - j}\right)^p$ |
| III | $\frac{e^{ja} n\beta^\gamma}{2\Gamma(\gamma)} \int_{-\infty}^{\infty} s^m e^{-(\beta s^n - jbs)} ds$ <br> $s = \frac{\ln y - a}{b}$ | $\gamma = \frac{m+1}{n} > 0,\ m > -1, \beta, n, b > 0, -\infty < a < \infty$ <br> EX: lognormal by GR 3.323.2, for $n = 2, m = 0, \beta = \frac{1}{2}$, <br> $\int_{-\infty}^{\infty} \exp(-p^2x^2 \pm qx)dx = \exp\left(\frac{q^2}{4p^2}\right)\frac{\sqrt{\pi}}{p}$; set $p^2 = \frac{1}{2}; q = j\sigma$, <br> giving $EY^j = \exp\left(j\mu + \frac{j^2\sigma^2}{2}\right)$ |
| IV | $\frac{e^{ja} n\beta^\gamma}{\Gamma(\gamma)} \int_{0}^{\infty} s^m e^{-(\beta s^n - jbs)} ds$ <br> $s = \frac{\ln y - a}{b}$ | $\gamma = \frac{m+1}{n} > 0,\ m > -1, \beta, n, b > 0, -\infty < a < \infty$ <br> EX: Pearson Type III mean for $n = 1, m = p - 1, \beta = 1$, <br> $\mu = \frac{e^a}{(1-b)^p}$ from $EY^j = \frac{e^{ja}}{(1 - jb)^p}$ |

- General solution is complicated and requires numerical methods.

EX. Derivation for Model II $EY^j$. By definition, $EY^j = \int_{-\infty}^{\infty} y^j g(y) dy$, so that, for $1 < y < \infty$,

$$EY^j = \int_{1}^{\infty} y^j \frac{n\beta^\gamma (\ln y)^m e^{-\beta(\ln y)^n}}{y\Gamma(\gamma)}, 1 < y < \infty.$$ Use change of variable: $\begin{cases} s = \ln y; ds = \frac{dy}{y}; s_1 = 0, s_\infty = \infty \\ e^s = y, y^j = e^{js} \end{cases}$

Substitute and obtain integral above which cannot be solved in closed form.

Example, Model III. The variance for the lognormal.

By definition, the variance is: $\sigma^2 = EY^2 - \mu^2$, where $\mu$ is the mean or $EY$ in above table, from derived expression, $EY^j = \exp\left(j\mu + \frac{j^2\sigma^2}{2}\right)$. The quantity $EY^2$ follows from $EY^j$ with $j$ =2. Expressed from the integral for the moments relation of the lognormal,

$$EY^2 = \frac{e^{2\mu}}{2\Gamma\left(\frac{1}{2}\right)} 2\left(\frac{1}{2}\right)^{\frac{1}{2}} \int_{-\infty}^{\infty} \exp-\left(\frac{s^2}{2} - 2\sigma s\right) ds = e^{2\mu+2\sigma^2}$$

Then the variance is: $\sigma^2 = e^{2\mu+2\sigma^2} - \left(e^{\mu+\frac{\sigma^2}{2}}\right)^2 = e^{2\mu+2\sigma^2} - e^{2\mu+\sigma^2}$. To obtain another form, simplify using one of the factoring rules for the difference of two exponentials,

$$e^x - e^y = \begin{cases} e^x\left(1 - e^{y-x}\right) \\ e^y\left(e^{x-y} - 1\right) \end{cases}$$

This gives another form for the lognormal variance:

$$\sigma^2 = e^{2\mu+2\sigma^2} - \left(e^{\mu+\frac{\sigma^2}{2}}\right)^2 = e^{2\mu+2\sigma^2} - e^{2\mu+\sigma^2} = \left(e^{\sigma^2} - 1\right)\left(e^{2\mu+\sigma^2}\right)$$

# Log CDFS

$F(y) = \int_{-\infty}^{y} g(t)dt$

| MODEL | Cumulative Distribution Function $F(y)$ | PARAMETER RESTRICTION | |
|---|---|---|---|
| I | $\frac{n\beta^\gamma}{2\Gamma(\gamma)}\int_0^y \frac{(\ln t)^m}{t} e^{-\beta(\ln t)^n} dt =$ $\frac{1}{2}\left[1+\frac{\gamma(\nu,\beta\ln y^n)}{\Gamma(\nu)}\right]$ | $\gamma = \nu = \frac{m+1}{n} > 0$ $m > -1, \beta, n > 0$ | • General solution for $F(y)$ requires series expansions on lower incomplete gamma function, $\gamma(a,x)$ <br>• Example: GEMII loggamma, $F(y) = \frac{\gamma(p,\lambda\ln y)}{\Gamma(p)}$ <br>• See Table 5a, NOTE 2, for limits assumptions for GEM III & IV CDFs <br>• See Table 5a, NOTE 7, for alternative expressions for CDFs <br>o E.g., GEM I, $F(y) = 1 - \frac{1}{2}\frac{\Gamma(\gamma,\beta\ln y^n)}{\Gamma(\gamma)}$ <br>• Complementary CDFs are: $1 - F(y) = P(Y > y) = \int_y^\infty f(t)dt.$ |
| II | $\frac{n\beta^\gamma}{\Gamma(\gamma)}\int_1^y \frac{(\ln t)^m}{t} e^{-\beta(\ln t)^n} dt =$ $\frac{\gamma(\nu,\beta\ln y^n)}{\Gamma(\nu)}$ | $\gamma = \nu = \frac{m+1}{n} > 0$ $m > -1, \beta, n > 0$ | |
| III | $\frac{n\beta^\gamma}{2b\Gamma(\gamma)}\int_0^y \left(\frac{\ln t - a}{b}\right)^m e^{-\beta\left(\frac{\ln t - a}{b}\right)^n} \frac{dt}{t} =$ $\frac{1}{2}\left\{1+\frac{\gamma\left[\nu,\beta\left(\frac{\ln y - a}{b}\right)^n\right]}{\Gamma(\gamma)}\right\}$ | $\gamma = \nu = \frac{m+1}{n} > 0$ $m > -1, \beta, b, n > 0,$ $-\infty < a < \infty$ | |
| IV | $\frac{n\beta^\gamma}{b\Gamma(\gamma)}\int_{e^a}^y \left(\frac{\ln t - a}{b}\right)^m e^{-\beta\left(\frac{\ln t - a}{b}\right)^n} \frac{dt}{t} =$ $\frac{\gamma\left[\nu,\beta\left(\frac{\ln y - a}{b}\right)^n\right]}{\Gamma(\nu)}$ | $\gamma = \nu = \frac{m+1}{n} > 0$ $m > -1, \beta, b, n > 0,$ $-\infty < a < \infty$ | |

- Example: derive the GEMIII lognormal CDF with parameters $[m = 0, n = 2, \beta = \gamma = 1/2]$;

$$F(y) = \frac{1}{\sigma\sqrt{2\pi}}\int_0^y e^{-\frac{1}{2}\left(\frac{\ln t - \mu}{\sigma}\right)^2} \frac{dt}{t} = \frac{1}{\sqrt{2\pi}}\int_{-\infty}^{\frac{\ln y - \mu}{\sigma}} e^{-\frac{1}{2}s^2} ds = 1 - \frac{1}{\sqrt{2\pi}}\int_{\frac{\ln y - \mu}{\sigma}}^{\infty} e^{-\frac{1}{2}s^2} ds = 1 - \frac{\Gamma\left[\frac{1}{2}, \frac{1}{2}\left(\frac{\ln y - \mu}{\sigma}\right)^2\right]}{2\sqrt{\pi}}$$

by change of variable, $s = \frac{\ln t - \mu}{\sigma}$, integration property of symmetric densities, $\int_{-\infty}^{y} f(t)dt = 1 - \int_y^\infty f(t)dt,$

and GR 3.381.9, $\int_u^\infty x^m e^{-\beta x^n} dx$

This leads to the error function,

$$F(y) = 1 - \frac{1}{2}\left[1 - \text{erf}\left[\frac{\ln y - \mu}{\sigma\sqrt{2}}\right]\right] = \frac{1}{2}\left[1 + \text{erf}\left(\frac{\ln y - \mu}{\sigma\sqrt{2}}\right)\right]$$

- GEMI lognormal $F(y)$ follows when $\mu = 0, \sigma = 1$.

Example: the median of the lognormal distribution (see discussion, Table 5a). By definition the median is the half way mark of the cumulative distribution. That is $F^{-1}\left(\frac{1}{2}\right) = \text{median}$. For the lognormal, this provides from the above result for the CDF, $F(y) = \frac{1}{2}\left[1 + \text{erf}\left(\frac{\ln y - \mu}{\sigma\sqrt{2}}\right)\right] = \frac{1}{2}$, and we must find $y$ in this relation that cuts the vertical axis $F(y)$ in half.

Algebraically simplifying the relation, $\text{erf}\left(\frac{\ln y - \mu}{\sigma\sqrt{2}}\right) = 0.$

It is known that the only value of the error function which gives $\text{erf}(x) = 0$ is $x = 0$.

That is, the definition of the error function (O'Brien, 2008) is:

$\text{erf}(x) = \frac{2}{\sqrt{\pi}}\int_0^x e^{-t^2}\,dt$. If $x = 0$, then $\text{erf}(0) = \frac{2}{\sqrt{\pi}}\int_0^0 e^{-t^2}\,dt = 0$. This means that $\frac{\ln y - \mu}{\sigma\sqrt{2}} = 0$ or $\ln y = \mu$, or $e^{\ln y} = e^{\mu} = y$ is the median of the lognormal. Thus, if $F(y) = \frac{1}{2}$, then $F^{-1}\left(\frac{1}{2}\right) = e^{\mu} = \text{median}$. Substitute $y = e^{\mu}$ into $F(y)$ to see that $\text{erf}(0) = 0$ and $F(y) = \frac{1}{2}$.

**Table A6.**

# Maximum Likelihood Estimation (MLE) Functions & Equations

| MODEL | PARAMETER SET $\theta$ | LIKELIHOOD FUNCTION $L(\Theta \mid x_1,\ldots,x_p)$ $p$ is sample size | LOG. LIKELIHOOD FUNCTION $\ln L(\Theta)$ $p$ is sample size |
|---|---|---|---|
| I<br>➢ Error<br>➢ Gaussian | $\gamma, n, \beta$ | $\prod_{i=1}^{p} \frac{n\beta^{\gamma}}{2\Gamma(\gamma)} x_i^{\,m} e^{-\beta x_i^{\,n}} = \left[\frac{n\beta^{\gamma}}{2\Gamma(\gamma)}\right]^{p} \prod_{i=1}^{p} x_i^{m} e^{-\beta\sum_{i=1}^{p} x_i^{n}}$ | $p\ln n + p\left(\frac{m+1}{n}\right)\ln\beta - p\ln 2 - p\ln\Gamma\left(\frac{m+1}{n}\right) + m\sum_{i=1}^{p}\ln x_i - \beta\sum_{i=1}^{p} x_i^{n}$ |
| | | | |
| II<br>➢ Exponential<br>➢ Gamma/ Erlang<br>➢ Ch-sq.<br>➢ Rayleigh<br>➢ Weibull<br>➢ Maxwell<br>➢ Nakagami<br>➢ Inverse gamma<br>➢ Generalized gamma<br>➢ Transformed gamma | $\gamma, n, \beta$ | $\prod_{i=1}^{p} \frac{n\beta^{\gamma}}{\Gamma(\gamma)} x_i^{\,m} e^{-\beta x_i^{\,n}} = \left[\frac{n\beta^{\gamma}}{\Gamma(\gamma)}\right]^{p} \prod_{i=1}^{p} x_i^{m} e^{-\beta\sum_{i=1}^{p} x_i^{n}}$ | $p\ln n + p\left(\frac{m+1}{n}\right)\ln\beta - p\ln\Gamma\left(\frac{m+1}{n}\right) + m\sum_{i=1}^{p}\ln x_i - \beta\sum_{i=1}^{p} x_i^{n}$ |
| | | | |
| III<br>➢ Gaussian<br>➢ Lognormal<br>➢ Laplace | $\gamma, n, \beta, a, b$ | $\prod_{i=1}^{p} \frac{n\beta^{\gamma}}{2b\Gamma(\gamma)} \left(\frac{x_i - a}{b}\right)^m e^{-\beta\left(\frac{x_i - a}{b}\right)^n}$<br>$= \left[\frac{n\beta^{\gamma}}{2b\Gamma(\gamma)}\right]^p \prod_{i=1}^{p} \left(\frac{x_i - a}{b}\right)^m e^{-\beta\sum_{i=1}^{p}\left(\frac{x_i - a}{b}\right)^n}$ | $p\ln n + p\left(\frac{m+1}{n}\right)\ln\beta - p\ln 2 - p\ln b -$<br>$p\ln\Gamma\left(\frac{m+1}{n}\right) + m\sum_{i=1}^{p}\ln\left(\frac{x_i - a}{b}\right) - \beta\sum_{i=1}^{p}\left(\frac{x_i - a}{b}\right)^n$ |
| | | | |
| IV<br>➢ Pearson (Type III)<br>➢ 2-parm. Exponential<br>➢ 3-parm. Weibul | $\gamma, n, \beta, a, b$ | $\prod_{i=1}^{p} \frac{n\beta^{\gamma}}{b\Gamma(\gamma)} \left(\frac{x_i - a}{b}\right)^m e^{-\beta\left(\frac{x_i - a}{b}\right)^n}$<br>$= \left[\frac{n\beta^{\gamma}}{b\Gamma(\gamma)}\right]^p \prod_{i=1}^{p} \left(\frac{x_i - a}{b}\right)^m e^{-\beta\sum_{i=1}^{p}\left(\frac{x_i - a}{b}\right)^n}$ | $p\ln n + p\left(\frac{m+1}{n}\right)\ln\beta - p\ln b -$<br>$p\ln\Gamma\left(\frac{m+1}{n}\right) + m\sum_{i=1}^{p}\ln\left(\frac{x_i - a}{b}\right) - \beta\sum_{i=1}^{p}\left(\frac{x_i - a}{b}\right)^n$ |

NOTES:

- MLE solution of Likelihood Equations, $\frac{\partial \ln \mathrm{L}(\Theta)}{\partial \Theta_i} = 0$, employs the transform $\frac{m+1}{n} \to \gamma$, and *m* is recovered from optimization solution of parameter set, $(\gamma, n, \beta)$.
- Odd/even functions algebra applies to Models I & III
- Suggested 3-parameter GEM II MLE solution (cf. Venter, 1983):

➢ Likelihood Equations, $\frac{\partial \ln \mathrm{L}(\gamma, n, \beta)}{\partial \Theta_i} = 0 \Rightarrow$

(a) $\frac{\partial \ln \mathrm{L}(\gamma, n, \beta)}{\partial \gamma} = p\ln\beta - p\Psi(\gamma) + n\sum_{i=1}^{p}\ln x_i = 0$

(b) $\frac{\partial \ln \mathrm{L}(\gamma, n, \beta)}{\partial n} = \frac{p}{n} + \gamma\sum_{i=1}^{p}\ln x_i - \beta\sum_{i=1}^{p} x^n \ln x = 0$

(c) $\frac{\partial \ln \mathrm{L}(\gamma, n, \beta)}{\partial \beta} = \frac{p\gamma}{\beta} - \sum_{i=1}^{p} x_i^n = 0$

- Solve $\beta$ in (c)
- Substitute $\beta$ into (b) and find $\gamma$
- Substitute $\beta$ and $\gamma$ into (a) to obtain one nonlinear equation in terms of *n*:

- $f(n)=\Psi(\gamma)-\ln(\gamma)-\frac{n}{p}\sum_{i=1}^{p}\ln x_i-\ln p+\ln\sum_{i=1}^{p}x_i^n=0$,

  where $\Psi(\bullet)$ is the Psi or Digamma Function, GR, Section 8.36, and

  - $\gamma=\dfrac{p\sum_{i=1}^{p}x_i^n}{n\left[p\sum_{i=1}^{p}x^n\ln x-\sum_{i=1}^{p}x_i^n\sum_{i=1}^{p}\ln x_i\right]}$

  - Solve for $n$ in one-variable Newton-Raphson optimization algorithm,
    - $n_{k+1}=n_k-\dfrac{f(n_k)}{f'(n_k)}, k=0,1,2,...$ with initial value, $n_0=2$.
      - The derivative is required for $f(n)$ where $f(n)$ w/o subscript $k$ on $n$ as follows:

$$f(n)=\Psi\left\{\frac{p\sum_{i=1}^{p}x_i^n}{n\left[p\sum_{i=1}^{p}x_i^n\ln x_i-\sum_{i=1}^{p}x_i^n\sum_{i=1}^{p}\ln x_i\right]}\right\}-\ln\left\{\frac{p\sum_{i=1}^{p}x_i^n}{n\left[p\sum_{i=1}^{p}x_i^n\ln x_i-\sum_{i=1}^{p}x_i^n\sum_{i=1}^{p}\ln x_i\right]}\right\}$$
$$-\frac{n}{p}\sum_{i=1}^{p}\ln x_i-\ln p+\ln\sum_{i=1}^{p}x_i^n=0$$

- GEM I MLE solution (symmetric density) is similar but odd-even algebra applies (see next §) to permissible solution set of $\gamma,n,\beta$.
  - ➢ NOTE: the literature has shown that this brute force approach does not work (i.e. does not converge because of parameter confounding). The log-exponential transforms (Table A5, NOTEs) are more amenable to nonlinear optimization solutions (see work on Stacy's 3-parm model in Lawless, 1980, 2003, Bell, 1988,and Noufaily & Jones, 2009), among others. Future work will take this approach for all models.

- The author concedes that Model II is very similar to Stacy 3-parameter model and Model IV is very similar to the Amoroso density if $\beta = 1$ which G. Crooks has shown subsumes many prob. distributions. The author's models were developed independent of knowledge of either but priority must be acknowledged.

- 5-parameter Models III & IV have not been solved but the log-exponential distributions (Notes, Table A5) discussed above are probably the preferred approach.

### ❖ Table A6a. Other Parameter Estimation Techniques

- Nonlinear regression
  - $\sum_{i=1}^{p} \left[ \ln y_i - \ln g(x_i) \right]^2 \to \min$

    where $g(x)$ is the modeling function of the Models I-IV.

◄TBD

- Minimun distance, minimum chi-square, percentile matching, methods of moments, graphical methods (Hogg & Klugman, 1984).

# Mathematics of Even-Odd Function PDFs

---

## A. Probability Density Function: Models I & II

Model I PDF: $f(x) = \frac{n\beta^{\gamma}}{2\Gamma(\gamma)} x^{m} e^{-\beta x^{n}}$ <> Model II PDF: $f(x) = \frac{n\beta^{\gamma}}{\Gamma(\gamma)} x^{m} e^{-\beta x^{n}}$

$-\infty < x < +\infty$ $\qquad\qquad$ $0 \le x < +\infty$

$\operatorname{Re} n > 0, \operatorname{Re}\beta > 0, \operatorname{Re}\gamma = \frac{m+1}{n} > 0, \operatorname{Re} m > -1$

- $m$ & $n$ have different meanings in each PDF IAW the values they can assume
- Model I PDF is an even function or symmetric density conceived as a generalized unimodal/bimodal Gaussian-like distribution on $\pm\infty$, with $\mu = 0, \sigma^2 > 0$.

Table A. Illustration of Model I (even function only) and Model II on Permissible PDF Parameter Values for *m* & *n* $(0 \le m, n \le 6)$ in Fractional Form

| $m = \frac{p}{q} > -1$ (reduced to lowest terms) | | | | $n = \frac{r}{s} > 0$ (reduced to lowest terms) | | | |
|---|---|---|---|---|---|---|---|
| $p$ | | $q$ | | $r$ | | $s$ | |
| Model I | Model II | Model I | Model II | Model I | Model II | Model I | Model II |
| 0 | 0 | | | | | | |
| | 1 | 1 | 1 | | 1 | 1 | 1 |
| 2 | 2 | | 2 | 2 | 2 | | 2 |
| | 3 | 3 | 3 | | 3 | 3 | 3 |
| 4 | 4 | | 4 | 4 | 4 | | 4 |
| | 5 | 5 | 5 | | 5 | 5 | 5 |
| 6 | 6 | | 6 | 6 | 6 | | 6 |

•NOTE: Definition of even function; $f(-x) = f(x)$, all $x$ (Bers, *Calculus*, p. 94).

- In probability theory, if a density is an even function, it is called a symmetric density based on a symmetric random variable $X$ (Hoel, et al., Vol. I, p. 123)

•NOTE: Model I PDF can be written: $f(x) \underset{n=\frac{r}{s}}{\overset{m=\frac{p}{q}}{=}} \frac{n\beta^{\gamma}}{2\Gamma(\gamma)} (x)^{\frac{p}{q}} e^{-\beta(x)^{\frac{r}{s}}}$

•NOTE: In each Model, negative values are allowed for $m = \frac{p}{q} > -1$ such as $-\frac{1}{3} - \frac{2}{3}, -\frac{4}{7}$, etc.

- NOTE: The $m = \frac{p}{q}$ parameter in Model I is restricted to the combinations of the set of all even integers (incl. 0) $\div$ all odd integers, and $n = \frac{r}{s}$ is the same (excluding 0). Model II has no such restrictions other than $m > -1, n > 0$.

Table B. Number of Logical Combinations of Permissible Values for *m* & *n* $(0 \le m, n \le 6)$

| Model | $m = \frac{p}{q} > -1$ | $n = \frac{r}{s} > 0$ |
|---|---|---|
| I | 12 | 9 |
| II | 42 | 36 |

NOTE: Numbers obtained from Table A by Fundamental Counting Principle (including duplicated values such as $\frac{2}{1}, \frac{6}{3}, \frac{4}{2}, \ldots$ ).
(Ignore fact that if *p*=0, *q* can be any even or odd integer > 0)

---

To illustrate the practical and theoretical differences between Model I & Model II, Table A shows the values that are permitted under Model I for parameters *m* & *n* if each cannot be greater than 6. Assume *m* is a fraction of integers $\frac{p}{q}$ (including 0) reduced to lowest terms , and *n* is a fraction of integers $\frac{r}{s}$ (excluding 0), reduced to lowest terms.

Table A show which values are allowed for *m* & *n* if the Model I PDF is to be an even function, or symmetric density, $f(-x) = f(x)$.

Table B indicates that if $(0 \le m, n \le 6)$, then the number of allowable combinations for Model II is much larger that of a Model I PDF—about 3-4 times as many for Model II. Formula to express difference_____? Thus, Model I is restricted in the permissible number of values that *m* & *n* can assume to achieve an even-function symmetric PDF.
Consequently,

- Model I PDF is a symmetric density and accepts only real-valued even functions on $\pm\infty$ with unimodal or bimodal graphs and mean of 0.
- Model II PDF accepts any real-valued function on $x \ge 0$ with unimodal graphs and mean $\neq 0$.

Similar arguments for Model III & IV as linear transforms of I/II

- Thus, modeling function $g(x) = x^3 e^{-\frac{1}{2}x^2}$ can be transformed into a Model II PDF but cannot be transformed into symmetric Model I PDF, because $x^3$ in $g(x)$ is not a symmetric random variable and the generated PDF will not be an even function or symmetric density.

---

## B. Technical Properties (Source: Bers, Hoel, et al., Wikipedia)

We note that an **even function** is symmetric about the $y$-axis, meaning that its graph remains unchanged after reflection about the $y$-axis; e.g., plots of Model I unimodal standardized Gaussian $N(0,1)$ PDF and GEMI$(0,\sigma^2)$ bimodal PDF, below,

Gaussian & Bimodal Even Function PDFs▶

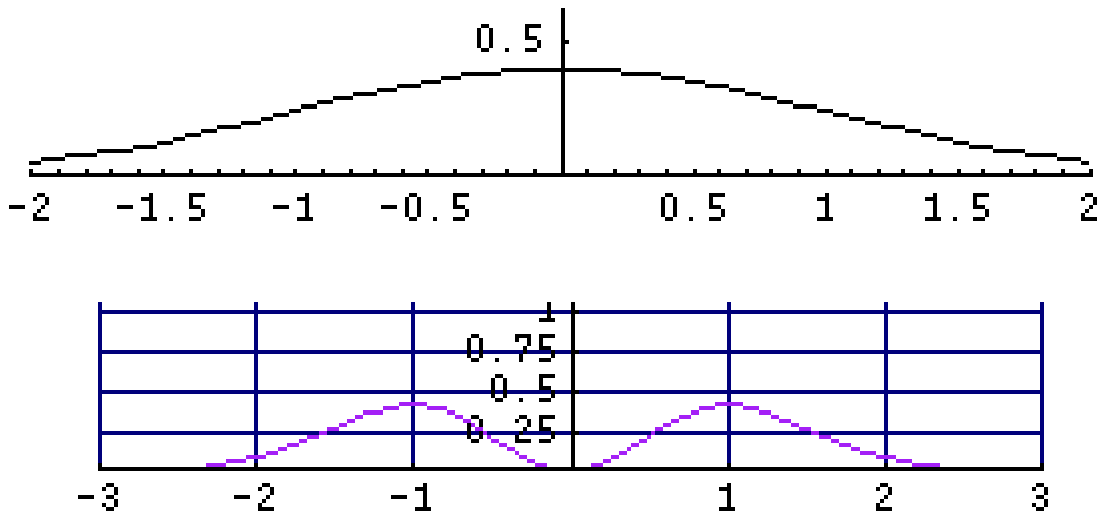

◀Bers, p. 94, other examples.

Fig. 1

**INTEGRAL PROPERTIES FOR SYMMETRIC PDF $f(x)$, GEMI ($\mu = \text{median} = 0$)**
**For GEMIII PDF replace** $0$ with $a$

$$f(-x) = f(x) \Rightarrow \int_{-\infty}^{\infty} f(x)dx = 2\int_{-\infty}^{0} f(x)dx = 2\int_{0}^{\infty} f(x)dx = 1$$

$$1 = \int_{-\infty}^{\infty} f(x)dx = \int_{-\infty}^{-x} f(x)dx + \int_{-x}^{0} f(x)dx + \int_{0}^{x} f(x)dx + \int_{x}^{\infty} f(x)dx = \int_{-\infty}^{x} f(x)dx + \int_{x}^{\infty} f(x)dx \Rightarrow \int_{-\infty}^{x} f(x)dx = 1 - \int_{x}^{\infty} f(x)dx$$

▲ Approach for GEMI cdf ▲

$$\int_{-\infty}^{0} f(x)dx = \int_{0}^{\infty} f(x)dx = \frac{1}{2} \Rightarrow \int_{-\infty}^{x} f(x)dx = \frac{1}{2} + \int_{0}^{x} f(x)dx$$

◀Showing second approach for $F(x)$

$$\int_{a}^{b} f(x)dx = \int_{-b}^{-a} f(x)dx, \text{ Bers, p. 261.}$$

$$\int_{-\infty}^{-x} f(x)dx = \int_{x}^{\infty} f(x)dx$$

◀Showing $F(-x) = 1 - F(x)$

Thus, CDFs—

$$F(x) = \int_{-\infty}^{x} = \left( \int_{-\infty}^{-x} + \int_{-x}^{0} + \int_{0}^{x} \right) = \left( \int_{-\infty}^{\infty} - \int_{x}^{\infty} \right) = 1 - \int_{x}^{\infty} = \frac{1}{2} + \int_{0}^{x}$$

$$F(-x) = \int_{-\infty}^{-x} = \int_{-\infty}^{\infty} - \left( \int_{-x}^{0} + \int_{0}^{x} + \int_{x}^{\infty} \right) = \int_{-\infty}^{\infty} - \int_{-x}^{\infty} = \int_{-\infty}^{\infty} - \int_{-\infty}^{x} = 1 - F(x)$$

$$\text{Check}: F(-x) + F(x) = \int_{-\infty}^{-x} + \int_{-\infty}^{x} = \int_{x}^{\infty} + \left( 1 - \int_{x}^{\infty} \right) = 1$$

An **odd function** is symmetric with respect to the origin, meaning that its graph remains unchanged after rotation of 180 degrees about the origin as in mean $(\mu = 0)$ of standardized Gaussian $N(0,1)$, below—

Odd Function ▶

$$\mu = \int_{-\infty}^{\infty} x f(x)\,dx = 0$$

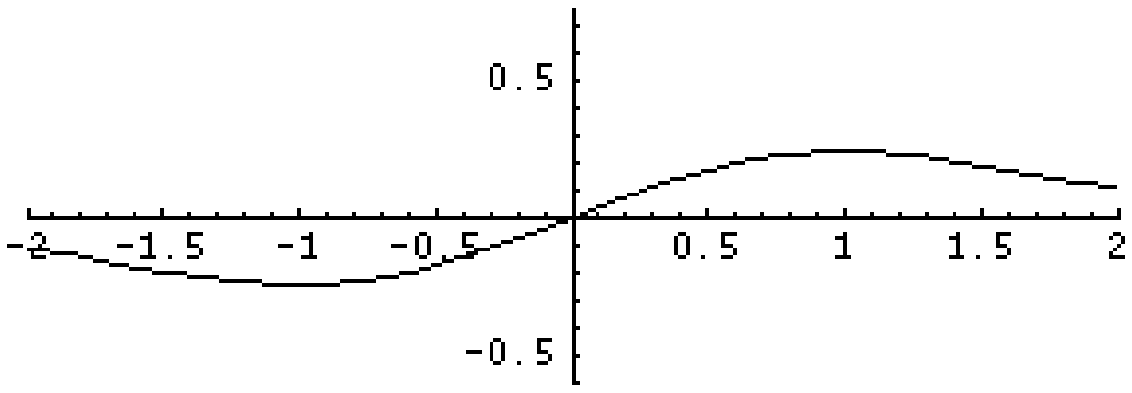

◀Bers, p. 94, other examples.

Fig. 2

**INTEGRAL PROPERTIES FOR ODD FUNCTIONS**

$$f(-x) = -f(x)$$

$$\int_a^b f(x)dx = -\int_{-b}^{-a} f(x)dx, \text{ Bers, p. 261.}$$

$$\int_{-a}^{a} f(x)dx = 0$$

$$\int_{-\infty}^{\infty} f(x)dx = 0$$

Compare following bimodal $(m \neq 0)$ & unimodal $(m = 0)$ Model I PDF graphs which are all even functions or symmetric densities (Purple curves are Gaussian), each with universal $\mu = 0$:

Model I
Even Function PDFs▶

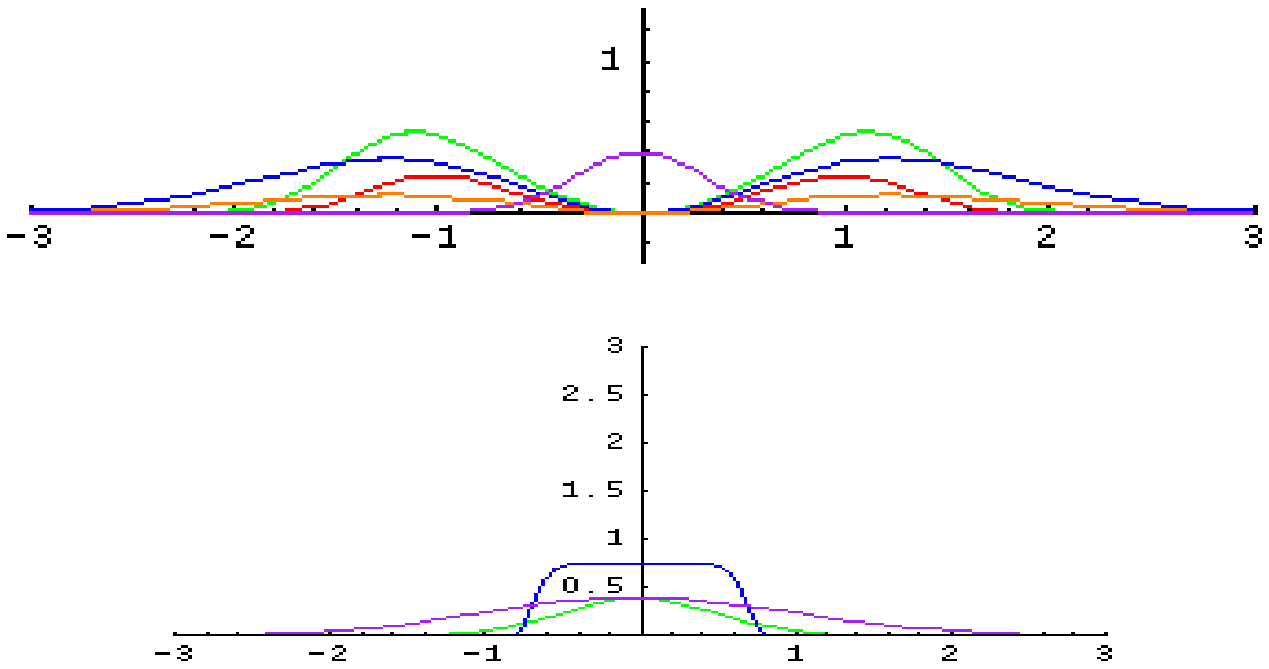

Fig. 3.

$f(x) = f(-x)$

▲Same principles apply for bimodal symmetric Even Function PDFs

$$\int_{-b}^{-a} f(x)dx = \int_a^b f(x)dx$$

$$\int_{-a}^{a} f(x)dx = 2\int_0^a f(x)dx$$

**C. Basic properties of even functions** (Source: Bers, Wikipedia)

- The only function which is *both* even and odd is the constant function which is identically zero (i.e., $f(x) = 0$ for all $x$).
- The sum of an even and odd function is neither even nor odd, unless one of the functions is identically zero.
- The sum of two even functions is even, and any constant multiple of an even function is even.
- The sum of two odd functions is odd, and any constant multiple of an odd function is odd.
- The product of two even functions is an even function.
- The product of two odd functions is again an even function.
- The product of an even function and an odd function is an odd function.
- The quotient of two even functions is an even function.
- The quotient of two odd functions is an even function.
- The quotient of an even function and an odd function is an odd function.

- The derivative of an even function is odd.
  - If $f(x) = \text{even}, f'(x) = \text{odd}; f''(x) = \text{even},\&\text{c}; f'(0) = -f'(0) = 0$ ◀Bers, pp. 544-5
- The derivative of an odd function is even.
  - If $f(x) = \text{odd}, f'(x) = \text{even}; f''(x) = \text{odd}, \&\ \text{c}, f'(0) = 0$ ◀Bers, pp. 544-5
  - Important for Mode calculation $\left(\frac{df(x)}{dx}\right)_{x=0}$.
- Th composition of two even functions is even, and the composition of two odd functions is odd.
- The composition of an even function and an odd function is even.
- The composition of any function with an even function is even (but not vice versa).
- The integral of an odd function from -A to +A is zero (where A is a finite, and the function has no vertical asymptotes between -A and A).
  - FOB: $\int_{-\infty}^{\infty} g(x)dx = 0$ if $g(x)$ is odd and converges—e.g., standardized Gaussian mean;Gradshteyn, I.S. and I.M. Ryzhik, 7th Edition, Sects. 3.03-3.05. ◀Fig. 2.
- The integral of an even function from -A to +A is twice the integral from 0 to +A (where A is a finite, and the function has no vertical asymptotes between -A and A).
  - FOB: $\int_{-\infty}^{\infty} g(x) = 2\int_{0}^{\infty} g(x)dx$ if $g(x)$ is even and converges—e.g., standardized Gaussian PDF;Gradshteyn, I.S. and I.M. Ryzhik, 7th Edition, Sects. 3.03-3.05 ◀Fig. 1.

## D. Series (Source: Bers & Wikipedia)

- The Taylor series of an even function includes only even powers.
  - "The Maclaurin series (Taylor series at $x = 0$) of an even function contains only even powers of $x$, whereas the Maclaurin series of an odd function contains only odd powers" (Bers, p. 545)
- The Taylor series of an odd function includes only odd powers.

- Taylor Series: $f(x)=f(x_0)+f'(x_0)(x-x_0)+\frac{f''(x_0)(x-x_0)^2}{2!}+\ldots$
  - Maclaurin: set $x_0 \to 0$.

◀Bers, pp. 544-5

- The Fourier series of a periodic even function includes only cosine terms.
- The Fourier series of a periodic odd function includes only sine terms.

## E. Algebraic structure (Source: Bers & Wikipedia)

- Any linear combination of even functions is even, and the even functions form a vector space over the reals. Similarly, any linear combination of odd functions is odd, and the odd functions also form a vector space over the reals. In fact, the vector space of *all* real-valued functions is the direct sum of the subspaces of even and odd functions. In other words, every function can be written uniquely as the sum of an even function and an odd function:
  - $f(x)=f_{\text{even}}(x)+f_{\text{odd}}(x)=\underbrace{\frac{f(x)+f(-x)}{2}}_{f_{even}}+\underbrace{\frac{f(x)-f(-x)}{2}}_{f_{odd}}=f(x)$

◀Example, Bers, p. 545

- The even functions form a commutative algebra over the reals. However, the odd functions do *not* form an algebra over the reals

## F. Model I Implications of Basic Properties: Programming Notes

- $EX^j=0$, $j$ odd
  - $EX^j=-EX^j=0$, $j$ odd
- $E(X-\mu)^j \overset{\mu=0}{=} 0$, $j$ odd
- $\sum_{k=1}^{p} x^k=0$, $k$ odd [$p$ is sample size]
- $\sum_{k=1}^{p} x^k \ln x=0$, $k$ odd [$p$ is sample size]
- $\ln x\sum_{k=1}^{p} x^k=0$, $k$ odd [$p$ is sample size]
- Other combinations such as quotients of odd/even simple functions $\left(\frac{x^3}{x^2}\right)$.
- Similar arguments for Models III.

Jointly Distributed Random Variables: Distributions of Sums and Quotients , Conditional and Marginal Densities, Covariance & Correlation, Multidimensional Change of Variable (Hoel, et al., Chap. 6 & Ch. 7, pp 178-180, 181ff.)

—TBD—

Hypothesis Testing. CLT Normal Approximation Formula, Kolmogorov-Smirnov, Other Nonparametric Tests, Chi-square

—TBD—

**TABLE A7**

# INTEGRAL FORMULAS

Source: *Table of Integrals, Series, and Products*. I.S. Gradshteyn and I.M. Ryzhik (7th Edition). Alan Jeffrey and Daniel Zwillinger, Editors. NY: Academic Press, 2007.

| GR FORMULA NUMBER[9] | FORMULA | CONSTRAINTS |
|---|---|---|
| 2.33.10 | $\int x^m e^{-\beta x^n} dx = -\frac{\Gamma(\gamma, \beta x^n)}{n\beta^\gamma}$ <br> $= -\frac{1}{n\beta^\gamma}\int_{\beta x^n}^{\infty} t^{\gamma-1}e^{-t}dt$ | $\gamma = \frac{m+1}{n}, n \neq 0, \beta \neq 0$ <br> For $\Gamma(\alpha, x)$ see 8.350.2 |
| — | $\int \left(\frac{x-a}{b}\right)^m e^{-\beta\left(\frac{x-a}{b}\right)^n} dx = -b\frac{\Gamma\left[\gamma, \beta\left(\frac{x-a}{b}\right)^n\right]}{n\beta^\gamma}$ <br> $= -\frac{b}{n\beta^\gamma}\int_{\beta\left(\frac{x-a}{b}\right)^n}^{\infty} t^{\gamma-1}e^{-t}dt$ | $\gamma = \frac{m+1}{n}, b \neq 0, n \neq 0, \beta \neq 0, \left\lvert\arg\left(\frac{x-a}{b}\right)\right\rvert > 0$ <br> For $\Gamma(\alpha, x)$ see 8.350.2 |
| 2.33.19 | $\int \frac{e^{-\beta x^n}}{x^m} dx = -\frac{\beta^z\Gamma(-z, \beta x^n)}{n}$ <br> $= -\frac{\beta^z}{n}\int_{\beta x^n}^{\infty} \frac{e^{-t}}{t^{z+1}}dt$ | $z = \frac{m-1}{n}, n \neq 0, \beta \neq 0$ <br> For $\Gamma(\alpha, x)$ see 8.350.2 |
| 3.323.2 | $\int_{-\infty}^{\infty} \exp\left(-p^2x^2 \pm qx\right)dx = \exp\left(\frac{q^2}{4p^2}\right)\frac{\sqrt{\pi}}{p}$ | $\operatorname{Re} p > 0$ |
| — | $\int_u^v \exp\left(\alpha x - \beta e^x\right)dx = \frac{1}{\beta^\alpha}\left[\Gamma(\alpha, \beta e^u) - \Gamma(\alpha, \beta e^v)\right]$ | Re $\beta > 0$ |

[9] Formulas identified by formula number in text (referred to as GR). A dash —in the first column indicates an unpublished formula. Known errata in textbook have been corrected in above table. Derivations in O'Brien (2008). Table includes long-standing GR formulas.

| | | |
|---|---|---|
| — | $\int_{-\infty}^{\infty} \exp(\mu x - \beta e^x) dx = \frac{1}{\beta^\mu}\Gamma(\mu)$ | $\mathrm{Re}\ \beta > 0$, $\mathrm{Re}\ \mu > 0$ |
| — | $\int_{-\infty}^{\infty} \exp(\alpha x - \alpha e^x) dx = \frac{1}{\alpha^\alpha}\Gamma(\alpha)$ | $\mathrm{Re}\ \alpha > 0$ |
| 3.326.2 | $\int_0^{\infty} x^m e^{-\beta x^n} dx = (n\beta^\gamma)^{-1}\Gamma(\gamma)$ | $\gamma = \frac{m+1}{n}$, $\mathrm{Re}\,\gamma > 0$, $\mathrm{Re}\,n > 0$, $\mathrm{Re}\,\beta > 0$ |
| 3.381.8 | $\int_0^{u} x^m e^{-\beta x^n} dx = (n\beta^\nu)^{-1}\gamma(\nu, \beta u^n)$ | $\nu = \frac{m+1}{n}$, $u > 0$, $\mathrm{Re}\,\nu > 0$, $\mathrm{Re}\,n > 0$, $\mathrm{Re}\,\beta > 0$ |
| 3.381.9 | $\int_u^{\infty} x^m e^{-\beta x^n} dx = (n\beta^\nu)^{-1}\Gamma(\nu, \beta u^n)$ | $\nu = \frac{m+1}{n}$, $\mathrm{Re}\,n > 0$, $\mathrm{Re}\,\beta > 0$ |
| 3.381.10 | $\int_0^{\infty} x^m e^{-\beta x^n} dx = \frac{\gamma(\nu, \beta u^n) + \Gamma(\nu, \beta u^n)}{n\beta^\nu}$ | $\nu = \frac{m+1}{n}$, $u > 0$, $\mathrm{Re}\,\nu > 0$, $\mathrm{Re}\,n > 0$, $\mathrm{Re}\,\beta > 0$<br>See also Formula **3.326.2** |
| — | $\int_u^{v} x^m \exp(-\beta x^n) dx = \frac{\Gamma(\gamma, \beta u^n) - \Gamma(\gamma, \beta v^n)}{n\beta^\gamma}$ | $\gamma = \frac{m+1}{n}$, $\mathrm{Re}\,\beta > 0$, $\mathrm{Re}\,n > 0$ |
| 3.381.11 (adapted) | $\int_{-\infty}^{\infty} x^m \exp(-\beta x^n) dx = \frac{2[\gamma(\nu, \beta x^n) + \Gamma(\nu, \beta x^n)]}{n\beta^\nu} = 2\frac{\Gamma(\nu)}{n\beta^\nu}$ | $\nu = \frac{m+1}{n}$, $\mathrm{Re}\,\beta > 0$, $\mathrm{Re}\ n > 0$, $\mathrm{Re}\ \nu > 0$<br>Even function |
| — | $\int_u^{v} \left(\frac{x-a}{b}\right)^m e^{-\beta\left(\frac{x-a}{b}\right)^n} dx = b\frac{\Gamma\left[\gamma, \beta\left(\frac{u-a}{b}\right)^n\right] - \Gamma\left[\gamma, \beta\left(\frac{v-a}{b}\right)^n\right]}{n\beta^\gamma}$ | $\gamma = \frac{m+1}{n}$, $-\infty < \mathrm{Re}\,a < +\infty$, $\mathrm{Re}\,b > 0$,<br>$\left\|\arg\left(\frac{u-a}{b}\right)\right\| > 0$, $\left\|\arg\left(\frac{v-a}{b}\right)\right\| > 0$, $\mathrm{Re}\,\beta > 0$, $\mathrm{Re}\,n > 0$ |
| 3.462.18 | $\int_0^{\infty} \left(\frac{x-a}{b}\right)^m e^{-\beta\left(\frac{x-a}{b}\right)^n} dx = \frac{b\Gamma\left[\gamma, \beta\left(-\frac{a}{b}\right)^n\right]}{n\beta^\gamma}$ | $\gamma = \frac{m+1}{n}$, $\mathrm{Re}\,b > 0$, $\arg(-a/b) > 0$, $\mathrm{Re}\,\beta > 0$, $\mathrm{Re}\,n > 0$ |
| — | $\int_{-\infty}^{\infty}\left(\frac{x-a}{b}\right)^m e^{-\beta\left(\frac{x-a}{b}\right)^n}dx=\frac{2b\Gamma(\gamma)}{n\beta^{\gamma}}$ | $\gamma=\frac{m+1}{n}$, $\operatorname{Re} b>0, \left\|\arg\left(\frac{x-a}{b}\right)\right\|>0, \operatorname{Re}\beta>0, \operatorname{Re} n>0, -\infty<a<\infty$ Even function |
| — | $\int_u^v \frac{1}{x^m}e^{-\beta x^n}dx=\beta^z\frac{\Gamma(-z,\beta u^n)-\Gamma(-z,\beta v^n)}{n}$ | $z=\frac{m-1}{n}, \operatorname{Re} -z>0, \operatorname{Re}\beta>0, \operatorname{Re} n>0$ |
| — | $\int_0^{\infty}\frac{e^{-\beta x^n}}{x^m}dx=\frac{\Gamma(z)}{n\beta^z}$ | $z=\frac{1-m}{n}, \operatorname{Re} z>0, \operatorname{Re}\beta>0, \operatorname{Re} n>0$ |
| 3.462.17 | $\int_u^{\infty}\frac{e^{-\beta x^n}}{x^m}dx=\frac{\Gamma(z,\beta u^n)}{n\beta^z}$ | $z=\frac{1-m}{n}, \operatorname{Re}\beta>0, \operatorname{Re} n>0$ |
| — | $\int_u^v x^m e^{-\frac{\beta}{x^n}}dx=\beta^{\gamma}\frac{\Gamma\left(-\gamma,\frac{\beta}{v^n}\right)-\Gamma\left(-\gamma,\frac{\beta}{u^n}\right)}{n}$ | $\gamma=\frac{m+1}{n}, \operatorname{Re}\beta>0, \operatorname{Re} n>0$ |
| — | $\int_0^{\infty}x^m e^{-\frac{\beta}{x^n}}dx=\beta^{\gamma}\frac{\Gamma(-\gamma)}{n}$ | $\gamma=\frac{m+1}{n}, \operatorname{Re} -\gamma>0, \operatorname{Re}\beta>0, \operatorname{Re} n>0$ |
| — | $\int_0^u\frac{e^{-\beta x^n}}{x^m}dx=\frac{\gamma(z,\beta u^n)}{n\beta^z}$ | $z=\frac{1-m}{n}, \operatorname{Re}\beta>0, \operatorname{Re} n>0$ |
| — | $\int_u^v\left(\frac{x-a}{b}\right)^m e^{-\frac{\beta}{\left(\frac{x-a}{b}\right)^n}}dx=b\frac{\Gamma\left(-\gamma,\frac{\beta}{\left(\frac{v-a}{b}\right)^n}\right)-\Gamma\left(-\gamma,\frac{\beta}{\left(\frac{u-a}{b}\right)^n}\right)}{n\beta^{-\gamma}}$ | $\gamma=\frac{m+1}{n}, -\infty<\operatorname{Re} a<\infty, \operatorname{Re} b>0, \left\|\arg\left(\frac{u-a}{b}\right)\right\|>0,$ $\left\|\arg\left(\frac{v-a}{b}\right)\right\|>0, \operatorname{Re}\beta>0, \operatorname{Re} n>0.$ |
| — | $\int_u^v\frac{e^{-\frac{\beta}{\left(\frac{x-a}{b}\right)^n}}}{\left(\frac{x-a}{b}\right)^m}dx=b\frac{\Gamma\left(z,\frac{\beta}{\left(\frac{v-a}{b}\right)^n}\right)-\Gamma\left(z,\frac{\beta}{\left(\frac{u-a}{b}\right)^n}\right)}{n\beta^z}$ | $z=\frac{m-1}{n}, -\infty<\operatorname{Re} a<\infty, \operatorname{Re} b>0, \left\|\arg\left(\frac{u-a}{b}\right)\right\|>0,$ $\left\|\arg\left(\frac{v-a}{b}\right)\right\|>0, \operatorname{Re}\beta>0, \operatorname{Re} n>0.$ |
| — | $\int_u^v \frac{e^{-\frac{\beta}{x^n}}}{x^m}dx = \frac{\Gamma\left(\frac{m-1}{n},\frac{\beta}{v^n}\right)-\Gamma\left(\frac{m-1}{n},\frac{\beta}{u^n}\right)}{n\beta^{\frac{m-1}{n}}}$ | $\operatorname{Re}\beta > 0, \operatorname{Re}n > 0$ |
| — | $\int_u^v \frac{e^{-\beta\left(\frac{x-a}{b}\right)^n}}{\left(\frac{x-a}{b}\right)^m}dx = b\frac{\Gamma\left(-z,\beta\left(\frac{u-a}{b}\right)^n\right)-\Gamma\left(-z,\beta\left(\frac{v-a}{b}\right)^n\right)}{n\beta^{-z}}$ | $z=\frac{m-1}{n}, -\infty < \operatorname{Re}a < \infty, \operatorname{Re}b > 0, \left\|\arg\left(\frac{u-a}{b}\right)\right\| > 0,$ $\left\|\arg\left(\frac{v-a}{b}\right)\right\| > 0, \operatorname{Re}\beta > 0, \operatorname{Re}n > 0.$ |
| — | $\int_0^\infty \frac{e^{-\frac{\beta}{x^n}}}{x^m}dx = \frac{\Gamma\left(\frac{m-1}{n}\right)}{n\beta^{\frac{m-1}{n}}}$ | $\operatorname{Re}\beta > 0, \operatorname{Re}n > 0, \operatorname{Re}\frac{m-1}{n} > 0$ |
| — | $\int_0^u \frac{e^{-\frac{\beta}{x^n}}}{x^m}dx = \frac{\Gamma\left(\frac{m-1}{n},\frac{\beta}{u^n}\right)}{n\beta^{\frac{m-1}{n}}}$ | $\operatorname{Re}\beta > 0, \operatorname{Re}n > 0$ |
| — | $\int_u^\infty \frac{e^{-\frac{\beta}{x^n}}}{x^m}dx = \frac{\Gamma\left(\frac{m-1}{n}\right)-\Gamma\left(\frac{m-1}{n},\frac{\beta}{u^n}\right)}{n\beta^{\frac{m-1}{n}}} = \frac{\gamma\left(\frac{m-1}{n},\frac{\beta}{u^n}\right)}{n\beta^{\frac{m-1}{n}}}$ | $\operatorname{Re}\beta > 0, \operatorname{Re}n > 0, \operatorname{Re}\frac{m-1}{n} > 0$ |
| 8.250.1<br>8.251.1 | $\frac{2}{\sqrt{\pi}}\int_0^x e^{-t^2}dt = \frac{1}{\sqrt{\pi}}\int_0^{x^2}\frac{e^{-t}}{\sqrt{t}}dt = \Phi(x)$ | $\Phi(x)$ is also called the error function, $\operatorname{erf}(x)$ |
| 8.310.1 | $\Gamma(z) = \int_0^\infty t^{z-1}e^{-t}dt$ | $\operatorname{Re}z > 0$ |
| 8.312.2 | $\Gamma(z) = x^z\int_0^\infty t^{z-1}e^{-xt}dt$ | $\operatorname{Re}z > 0, \operatorname{Re}x > 0$ |
| — | $\Gamma(-z) = -\frac{\Gamma(1-z)}{z}$ | $\operatorname{Re}z \neq 0,-1,-2,-3,\ldots$ |
| 8.350.1 | $\gamma(\alpha,x) = \int_0^x t^{\alpha-1}e^{-t}dt$ | $\operatorname{Re}\alpha > 0$ |
| — | $\gamma(\alpha, x) = x^{a}\int_{0}^{1} t^{\alpha-1}e^{-xt}dt$ | $\mathrm{Re}\,\alpha > 0, \mathrm{Re}\,x > 0$ |
| 8.350.2 | $\Gamma(\alpha, x) = \int_{x}^{\infty} t^{\alpha-1}e^{-t}dt$ | |
| — | $\Gamma(\alpha, x) = x^{a}\int_{1}^{\infty} t^{\alpha-1}e^{-xt}dt$ | $\mathrm{Re}\,\alpha > 0, \mathrm{Re}\,x > 0$ |
| 8.350.3 | $\Gamma(\alpha,0) = \Gamma(\alpha)$ | |
| 8.350.4 | $\Gamma(\alpha,\infty) = 0$ | |
| 8.350.5 | $\gamma(\alpha,0) = 0$ | |
| — | $\gamma(\alpha,\infty) = \Gamma(\alpha)$ | |
| 8.352.6 | $\gamma(n, x) = (n-1)!\left[1 - e^{-x}\left(\sum_{m=0}^{n-1}\frac{x^m}{m!}\right)\right]$ | $n = 1,2,...$ |
| 8.352.7 | $\Gamma(n, x) = (n-1)!\, e^{-x}\sum_{m=0}^{n-1}\frac{x^m}{m!}$ | $n = 1,2,...$ |
| 8.354.1 | $\gamma(\alpha, x) = \sum_{n=0}^{\infty}\frac{(-1)^n x^{\alpha+n}}{n!(\alpha+n)}$ | |
| 8.354.2 | $\Gamma(\alpha, x) = \Gamma(\alpha) - \sum_{n=0}^{\infty}\frac{(-1)^n x^{\alpha+n}}{n!(\alpha+n)}$ | |
| 8.354.5 | $\Gamma(\alpha, x) = e^{-x}x^{\alpha}\sum_{n=0}^{\infty}\frac{L_n^{\alpha}(x)}{n+1}$ | $x > 0$<br>$L_n^{\alpha}(x)$ are Laguerre polynomials, Sect. 8.97 |
| 8.356.3 | $\Gamma(\alpha) = \gamma(\alpha, x) + \Gamma(\alpha, x)$ | |
| 8.356.4 | $\frac{d\gamma(\alpha, x)}{dx} = -\frac{d\Gamma(\alpha, x)}{dx} = x^{\alpha-1}e^{-x}$<br>$\frac{d\gamma(a,u)}{dx} = -\frac{d\Gamma(a,u)}{dx} = \frac{d\gamma(a,u)}{du}\frac{du}{dx}$ | |
| 8.359.3 | $\Gamma\left(\frac{1}{2}, x^2\right) = \sqrt{\pi} - \sqrt{\pi}\Phi(x)$ | $\Phi(x)$ is probability integral (or error function), Form. 8.250.1 |
| 8.359.4 | $\gamma\left(\frac{1}{2}, x^2\right) = \sqrt{\pi}\Phi(x)$ | $\Phi(x)$ is probability integral (or error function), Form. 8.250.1 |
| — | $\Gamma\left(\frac{1}{2}, x\right) = \sqrt{\pi} - \sqrt{\pi}\Phi\left(\sqrt{x}\right)$ | $\Phi(x)$ is probability integral (or error function), Form. 8.250.1 |
| — | $\gamma\left(\frac{1}{2}, x\right) = \sqrt{\pi}\,\Phi\left(\sqrt{x}\right)$ | $\Phi(x)$ is probability integral (or error function), Form. 8.250.1 |
| — | $\Gamma\left(\frac{n}{2}, x\right) = \frac{(n-2)!!}{2^{\frac{n-3}{2}}}\left[\frac{1}{2}\Gamma\left(\frac{1}{2}, x\right) + e^{-x}x^{\frac{1}{2}}\sum_{s=0}^{\frac{n-3}{2}}\frac{(2x)^s}{(2s+1)(2s-1)!!}\right]$ $[n = 1,3,5,\ldots]$ | |
| — | $\gamma\left(\frac{n}{2}, x\right) = \frac{(n-2)!!}{2^{\frac{n-3}{2}}}\left[\frac{1}{2}\sqrt{\pi}\Phi\left(\sqrt{x}\right) - e^{-x}x^{\frac{1}{2}}\sum_{s=0}^{\frac{n-3}{2}}\frac{(2x)^s}{(2s+1)(2s-1)!!}\right]$ $[n = 1,3,5,\ldots]$ | |

# # #

# NOTES